\documentclass[11pt]{article}
\usepackage{amsfonts}
\usepackage{amsthm}
\usepackage{amssymb}
\usepackage{amsmath,amsthm}
\usepackage{amsthm}
\usepackage{enumerate}
\usepackage{amssymb}
\usepackage{fancyhdr}
\usepackage[dvips]{graphicx}
\newtheorem{theorem}{Theorem}[section]
\newtheorem{proposition}[theorem]{Proposition}
\newtheorem{corollary}[theorem]{Corollary}
\theoremstyle{definition}
\newtheorem{definition}[theorem]{Definition}
\newtheorem{assumption}[theorem]{Assumption}
\newtheorem{lemma}[theorem]{Lemma}
\newtheorem{example}[theorem]{Example}
\newtheorem{remark}[theorem]{Remark}
\newtheorem{hypothesis}[theorem]{Hypothesis}

\begin{document}
\title{How do future expectations affect the financial sector? Expectations modeling, infinite horizon (controlled) random FBSDEs and stochastic viscosity solutions.\footnote{The authors wish to thank Professor Josef Teichmann for his fruitful suggestions. This work was supported by Research Program THALES ``Support of the interdisciplinary and/or inter-institutional research and innovation", under the auspices of National Strategic Reference Framework (NSRF) 2007-13.}}
\author{Xanthi-Isidora Kartala$^{(a)}$, Nikolaos Englezos$^{(b)}$ and Athanasios N. Yannacopoulos$^{(a)}$   \\  \\ (a) Department of Statistics, Athens University of Economics and Business. \\ (b) Department of Banking and Financial Management, University of Piraeus.}
\date{}
\maketitle

\begin{abstract}
In this paper we study a class of infinite horizon fully coupled forward-backward stochastic differential equations (FBSDEs), that are stimulated by various continuous time future expectations models with random coefficients.
Under standard Lipschitz and monotonicity conditions, and by means of the contraction mapping principle, we establish existence, uniqueness, a comparison property and dependence on a parameter of adapted solutions.
Making further the connection with infinite horizon quasilinear backward stochastic partial differential equations (BSPDEs) via a generalization of the well known four-step-scheme,
we are led to the notion of stationary stochastic viscosity solutions. A stochastic maximum principle for the optimal control problem of such FBSDEs is also provided as an application to this framework.
\end{abstract}
\medskip

\medskip


\section{Introduction}\label{sec:intro}

Expectations concerning the future state of the economy are long  known to affect its present state as well. Leading economists have tried to model and understand the effect of future expectations to today's economy as well as their long run effects. A key role in this theory has been played by rational expectations theory, according to which the economic agents act based on their expectations which are
updated in a ration fashion to match actual facts, and this repeated activity is assumed to drive the system to some equilibrium. However, future expectations can be conceived as having the opposite effect, i.e., destabilizing the system. For example, one may argue that the wrong expectations (e.g., over-optimistic expectations concerning the market) may lead to instability and bubble formation.
These arguments have led to a number of random dynamical models, of varying mathematical sophistication, whose key idea is that the agents observe signals according to which they try to predict future events in the economy and base their actions on these predictions. The degree to which this signal may carry relevant information for the prediction of the quantities of interest or may be interpreted by agents in a ration fashion to be of some use can vary.

The economy is a live, complex and evolving system, comprising of agents whose actions are affected but at the same time affect the state of the system. The main idea of this paper is that there are quantities (often called the \emph{fundamentals} of the economy) whose evolution affects the actions of the agents, who often act in response to their predictions of the future values of such quantities. The present actions of the agents affect a sector of the economy (in our case the financial sector) and may cause fluctuations of the prices of the  \emph{assets} in this sector, which in turn affect the actual (as opposed to the predicted by the agents) value of the fundamentals. Then, assuming agents to be rational, their predictions for future values of the fundamentals are re-evaluated leading to a change in their course of action and so on and so forth. The way a control engineer would pose it is that the economy is considered as a feedback loop, "designed" by the invisible hand so as to equilibrate in the long run as a result of this scheme, which models the interplay of future expectations with today's actions. The way a mathematician would pose it would be to consider it as a dynamical system (possibly random) which has a fixed point (hopefully stable) that is a desired state of the economy. The way an economist would pose it, agents' individual actions will in the long run equilibrate the economy (however, as another economist may state it "In the long run we are all dead.").

By persuasion and training we shall take up the second stand and try to model the above situation as a \emph{random dynamical system}. The state space would consist of two contributions $(F,A)$ where $F$ are the values of the  fundamentals (i.e. the signal) of the economy and $A$ are the value of the assets. Their evolution will be given in a schematic fashion as  $F_{p}={\cal G}_{1}(A_{p},\omega)$ and $A_{p}={\cal G}_{2}(F_{p},_{f},\omega),$ where by the subscript $p$ we denote the value of the relevant quantity today (at present) whereas by the subscript $f$ we denote a prediction (an expectation) of the future values of this quantity based on the information available at present. The exact forms of the functionals ${\cal G}_{i}$, $i=1,2,$ are to be determined, however it is important to  note that they can be random (denoted by the inclusion of $\omega$); a fact that may be used to model a variety of effects such as departures from rationality of the agents, the basement of predictions on statistical or econometric tests etc. A fixed point of this dynamical system, if it exists, would be the desired \emph{equilibrium}.
Another desired output would be to obtain a relation of the form $A_{p}={\cal G}(F_{\le p}),$ where $F_{\le p}$ denotes the past and the present (but not the future)  values of the fundamentals, which will act as a rational \emph{pricing rule}.

Since our motivation is finance, we start our analysis by mentioning cases where such a model can be useful, e.g., in studying  the interaction of the macro economy with the financial sector (where $F$ is
a macroeconomic factor and $A$ is some index of the financial market), in studying  exchange rate dynamics (where $F$ is related to some policy variable as the domestic price level and $A$ is the foreign exchange rate)
or in studying the effects of pricing, subject to future expectations concerning the market, may have on the stabilization or destabilization of the market. In particular, we make a short review
of several models of this kind, such as the Krugman's \cite{Krugman} target zones exchange rate model, the Dornbusch's \cite{Dornbusch} model of exchange rate overshooting, the Blanchard's \cite{Blanchard} model for the interaction of the real economy with the stock market, and the Black's consol rate conjecture by Duffie et al \cite{Duffie-Ma-Yong} for the interplay between the short rate and perpetual bonds.

All the above models assumed that the coefficients of the saddlepoint dynamical systems were not random. However, since these models may involve predictions obtained by time series or econometric tools, or even possible
deviations of the agents from rationality, it seems reasonable to assume that this situation can be modeled by coefficients which are no longer deterministic functions of the variables but rather \emph{random fields}.
Such a generalization will  be important in order to introduce more realistic effects and produce qualitative results towards convergence to an equilibrium which are more robust and less sensitive to the actual assumptions
of the model.
% While this introduces more realistic effects in the model which may well affect its qualitative behavior (such as stability or existence of equilibrium), at the same time it complicates it by introducing a number of difficulties,
% among which the failure of being able to use the four-step-scheme in a form that provides a deterministic elliptic PDE.
Hence, we present a general future expectations model that covers all the above cases and show that, under simple and natural monotonicity and regularity conditions on its random coefficients, it is equivalently formulated
as an infinite time-horizon random forward-backward stochastic differential equations (FBSDEs) system. Then it is clear that the \emph{well-posedness} of the future expectations model, and thus of all the above, boils down to the solvability and the qualitative properties of the
associated system of random FBSDEs.

In this vein, we establish a more general form of an infinite horizon system of FBSDEs, whose random coefficients satisfy now a generalized version of the aforementioned conditions. In order for this system not to explode but reach eventually equilibrium after long enough time, we also impose a square-integrability condition to its solution with respect to a scalar discount rate of a given range. Then, making use of the contraction mapping principle, we extend Yin's \cite{Yin} result (cf. Theorem 3.1 therein) on \emph{existence} and \emph{uniqueness} of solutions of such FBSDEs by providing a more general lower bound for the above range. Another result of this analysis is an explicit
relationship that shows exactly how the stochastic discount factor of the primary future expectations model is confined by this scalar discount rate. In addition, under the imposed monotonicity and regularity conditions to
the random coefficients of these infinite horizon FBSDEs, we obtain a comparison result in case of an one-dimensional backward process at initial time, as well as continuous dependence of the solution on its initial conditions, i.e., the initial time and state of the forward process. Together with the solvability of the corresponding system of infinite horizon FBSDEs, these qualitative properties ensure the well-posedness of the future expectations model as they establish the price sensitivity of the assets in the economy, represented by the backward process, in terms of the fundamentals which are formulated by the forward process of the saddlepoint system.

Furthermore, we are interested in finding the \emph{representing field} that expresses instantaneously the backward variable of the above FBSDEs system in terms of the forward one; then the system would be decoupled, since its forward equation would reduce to a forward stochastic differential equation (FSDE) with respect only to the forward process, whose solvability is now easier to treat, and the backward processes would follow from the representing field. Considering the finite time-horizon analogue of the system, we show via the generalized It\^{o}-Kunita-Wentzell (IKW) formula and the well known four-step-scheme proposed by Ma and Yong \cite{Ma-Yong98} that this representing field is actually the unique random field pair which solves, in the classical sense, an associated  quasilinear second-order parabolic backward stochastic partial differential equation (BSPDE). However, since it is not always possible to obtain sufficient regularity for a random field pair to be a classical solution of this BSPDE, we are led next to a generalized definition of a \emph{stochastic viscosity solution} which to the best of our knowledge is new to the existing literature. In fact, it turns out that the backward process of the FBSDEs system yields such a solution to the BSPDE. If this solution is also unique then it provides the unknown initial condition of the backward process in terms of the known initial condition of the forward one, allowing us to treat the FBSDEs system only as an FSDEs system that
could be solved via a shooting algorithm. Getting eventually back to the infinite horizon case, we establish the conditions that should be imposed on the coefficients of the FBSDEs system under which the above stochastic
viscosity solution becomes \emph{stationary} and accordingly has the same probability distribution at all times.

Finally, we may revisit the aforementioned exchange rate economic models of Krugman \cite{Krugman} and Dornbusch \cite{Dornbusch}, where now an exogenous market variable as the domestic money supply is assumed to be determined by the economic authorities, in order for the resulting exchange rate to achieve a certain target subject to a costly intervention mechanism. Given that these economic models gave raise to future expectations models and in turn to infinite horizon systems of random FBSDEs, whose solvability has already been established, we are naturally led to the investigation of the \emph{stochastic control} problem for these systems. Hence, we consider a \emph{state} infinite horizon FBSDEs system whose random coefficients depend on a control process that takes values in a general domain and must be chosen such as the associated exponentially discounted
 \emph{cost functional} is minimized. By means of the duality theory, we introduce the generalized \emph{Hamiltonian random field} in order to state the \emph{adjoint} infinite horizon system of FBSDEs, according to which
 we obtain eventually a \emph{stochastic maximum principle} for the existence of an \emph{optimal control} process. 
\smallskip

\emph{Preview}. In Section \ref{sec:rem}, we introduce random effects as random fields in several well known future expectation models from economics or finance and express them in terms of infinite horizon FBSDEs systems with random coefficients. Then, through this reformulation, we study in Section \ref{sec:ssa} the long term behavior of these systems by providing solvability results in a proper functional setting for the resulting random FBSDEs. A representation of the functional relation $A_{p}={\cal G}(F_{\le p})$ is provided in Section \ref{sub:svs} by showing that now ${\cal G}$ is a random field which is the solution of a BSPDE of parabolic form. This also leads to the introduction of a concept of stable solutions for such BSPDEs which are stochastic generalizations of viscosity solutions. Finally, Section \ref{sec:ihmp} considers the case where the controlled
system of FBSDEs can reach a desired state by the minimum amount of intervention, giving rise to the development of an appropriate stochastic maximum principle.

\smallskip

\emph{Literature review.} Economic and financial models of future expectations have been proposed by a number of authors,  in a variety of fields and via a wide range of mathematical sophistication. Among some very influential papers in this direction   we can mention Dornbusch \cite{Dornbusch},  Blanchard \cite{Blanchard},  Krugman \cite{Krugman} and   Miller and  Weller \cite{Miller-Weller}, who introduced forward looking dynamical systems in trying to model effects of future expectations in various fields related to finance or macro-economics. From the stochastic analysis point of view, the convenient framework to study this type of forward looking dynamics is when randomness is introduced in the model by stochastic factors of semimartingale type.  On a related strand,  of great importance has been the works of  Duffie et al \cite{Duffie-Ma-Yong} who  considered the celebrated Black's consol rate conjecture as a forward looking system, which was importantly expressed in the form of a backward stochastic differential equation (BSDE), and Cvitani\'{c} and  Ma \cite{Cvitanic-Ma} who studied the effects of future expectations of large investors in the market in terms of FBSDEs. Furthermore, Yannacopoulos \cite{Yannacopoulos} reformulated the general class of models in Miller and Weller \cite{Miller-Weller}  as a system of FBSDEs and provided general results regarding their convergence towards equilibrium. In these works the celebrated four-step-scheme for the study of FBSDEs (cf. \cite{Ma-Yong98}) has played an important role in obtaining the functional relation $A_{p}={\cal G}(F_{\le p})$ by identifying it as the solution of a deterministic elliptic partial differential equation (PDE), thus providing important information on stability and instability (bubble formation) effects or possible pricing rules.

Existence and uniqueness results of an adapted solution for random systems of FBSDEs, under a traditional monotonicity condition on
their coefficients, were given via a \emph{homotopy method} by
Peng and Wu \cite{Peng-Wu}
 %\cite{Peng-Wu}
with a finite time-horizon and by Peng and Shi \cite{Peng-Shi}
in the infinite horizon case. Under a different but just as simple and natural monotonicity condition,
 Pardoux and Tang \cite{Pardoux-Tang}
and Yin \cite{Yin}
proved by means of the \emph{contraction mapping principle} the solvability of such systems in the finite and infinite time-horizon case, respectively, as well.

In a Markovian framework, Cvitani\'{c} and Ma \cite{Cvitanic-Ma}
 proved a comparison theorem for adapted solutions of FBSDEs systems
via the four-step-scheme and PDE methods. In the non Markovian case, comparison results for finite horizon FBSDEs
were obtained by Wu \cite{Wu99}
at initial time when the backward process is one-dimensional, and by Wu and Xu \cite{Wu-Xu}
at any time in higher dimensions. As regards an infinite time-horizon, Peng and Shi \cite{Peng-Shi}
gave comparison theorems for all times in the one dimensional case, while recently Zhang and Shi \cite{Zhang-Shi}
extended these results in the multidimensional one.

The notion of viscosity solutions for finite horizon deterministic PDEs was initially proposed by Crandall and Lions \cite{Crandall-Lions}
to amend the lack of regularity displayed by the value function of an optimal stochastic control problem in order to constitute a solution of
the associated Hamilton-Jacobi-Bellman (HJB) PDE. For a detailed exposition of this theory we refer to the "User's guide" paper by Crandall et al. \cite{Crandalletall}
and the monograph by Fleming and Soner \cite{Fleming-Soner}.
Lions and Souganidis \cite{Lions-Souganidis}, \cite{Lions-Souganidis98}
extended for the first time the notion of viscosity solutions to stochastic partial differential equations (SPDEs); in particular, they defined a weak solution for an SPDE
that is satisfied by the uniform limit of all its possible deterministic approximations. In the same spirit,  Buckdahn and Ma \cite{Buckdahn-Ma}, \cite{Buckdahn-Ma01},
defined a stochastic viscosity solution for an SPDE, which via a non-linear Doss-Sussmann-type transformation can be converted
to an ordinary PDE with random coefficients and thus can be studied $\omega$-wisely as a deterministic one. Employing these results,
Buckdahn and Ma \cite{Buckdahn-Ma07}
characterized the value function of a class of pathwise stochastic control problems as the unique stochastic viscosity solution, in their sense,
of the arising stochastic HJB equation.

The work of Dupire \cite{Dupire}
 on a new type of stochastic calculus for non-anticipative path functionals, which was further elaborated by
 Cont and Fournie \cite{Conte-Fournie}, \cite{Conte-Fournie10}
with applications to finance, has triggered recent developments in the generalization of the viscosity notion to path-dependent PDEs (PPDEs). Introducing a "left frozen maximization" approach for path functionals that makes the space of paths compact, Peng  \cite{Peng11}
obtained the maximum principle of fully non-linear PPDEs for a new notion of a stochastic viscosity solution. Alternatively, in two chain projects Ekren et al. \cite{Ekren-Keller-Touzi-Zhang}, \cite{Ekren-Touzi-Zang}, \cite{Ekren-Touzi-Zang12} proved uniqueness, existence, stability and comparison principle of another notion of stochastic viscosity solutions for either a semi-linear or a fully non-linear PPDE, respectively, on the space of continuous paths by Peron's method without compactness arguments.
Restricting the semi-jets on a family of expanding compact subsets of an $\alpha$-Holder space $C^\alpha$ for $\alpha\in(0,\frac{1}{2})$, a novel definition of stochastic viscosity solutions was also given in Tang and Zhang \cite{Tang-Zang}
in order to prove the verification theorem for the resulting fully non-linear Bellman PPDE of a path dependent optimal stochastic control problem.

In the case of a finite time-horizon, a stochastic hamilton system of FBSDEs arises in Peng and Wu \cite{Peng-Wu}
for a stochastic optimal control problem of an FSDE, subject to the aforementioned traditional monotonicity condition. Under the latter
a \emph{necessary} condition, known as the \emph{Pontryagin's} \emph{maximum principle}, for an optimal control of fully coupled FBSDEs was provided
by variational and duality techniques in  Wu \cite{Wu}
given a convex control domain and, alternatively, in Shi and Wu \cite{Shi-Wu}
given a control domain that is not necessarily convex but the diffusion coefficient does not depend on the control variable.  A unifying result
can be found in Yong \cite{Yong}
for controlled FBSDEs with randomly disturbed mixed-terminal conditions by means of spike variation techniques.

 In the infinite
time-horizon case, maximum principles were established for the stochastic control problem of a jump diffusion with partial information
by Haadem et al. \cite{Haadem-Oksendal-Proske}
and of an one dimensional FBSDEs system with jumps and delay by Agram and {\O}ksendal \cite{Agram-Oksendal},
subject to appropriate \emph{transversality conditions}.
Under the different monotonicity condition mentioned above, a \emph{sufficient} version of the maximum principle was given in Maslowski and Veverka \cite{Veverka}
for the discounting control problem of an FSDE by a proper approximation
of the backward adjoint process.

\section{Future Expectations Models And Infinite Horizon Random FBSDEs: Mathematical Preliminaries and Motivating Examples}\label{sec:rem}% And Review Of Related Literature

We start this section by introducing the general mathematical framework of the paper. Then we provide motivation for the connection of future expectations models with systems of FBSDEs by revisiting some well known and established in economics and finance models of this form and arguing as
to the benefit of allowing the coefficients to be random; a fact that as we shall see leads us naturally to the study of infinite horizon random FBSDEs.
%To set notation, let throughout this chapter $(\Omega ,\mathcal{F}, P; \mathbb{F})$ be a filtered
%probability space with $\mathbb{F}\triangleq\big(\mathcal{F}\big(t)\big)_{t\geq 0}$ being the natural filtration, generated by a $d-$dimensional Wiener process  $W(\cdot),$
%augmented by all the $P-$ null sets in $\mathcal{F}.$

%\subsection{Motivations From Mathematical Economics And Finance}\label{sub:MMEF}  7777
\vskip10pt

\subsection{Mathematical Preliminaries}\label{subsec:PRE}

Throughout this paper,
let $\big(\Omega, \mathcal{F},P; \mathbb{F}\big)$
be a filtered probability space carrying a standard $d$-dimensional Brownian motion $W(\cdot)$ and let $\mathbb{F}\triangleq\big\{\mathcal{F}(t); t\geq 0\big\}$ be the $\sigma$-field generated by $W(\cdot),$ that is $\mathcal{F}(t)=\sigma\{W(s): 0\leq s \leq t\}.$ We make also the standard $P$-augmentation to each $\mathcal{F}(t)$ such that it contains all the  $P$-null sets of $\mathcal{F}.$

Let $\lessdot,\gtrdot ,$  $|\cdot|$ and $||\cdot||$ stand for the usual scalar product, the Euclidean norm for vectors and the trace norm for matrices, respectively. Furthermore, let
 \verb"x"$^{\top}$ denote the transpose of a vector or matrix \verb"x", $\nabla_.$ stand for the gradient operator, subscripts denote partial derivatives with respect to the corresponding variable,
 and $a \vee b$ denote the maximum of $a$ and $b$. Moreover, comparison of vectors should be understood pointwisely between their respective coordinates.

For any integer $k\geq 0$ and any given $T>0$, we denote by $C^{0,k}\big([0,T]\times\mathbb{R}^{n};\mathbb{B}\big)$
 the set of all $\mathbb{B}$-valued functions
on $[0,T]\times\mathbb{R}^{n}$ that are continuous with respect to $t \in [0, T]$
and continuously
differentiable up to order $k$ with respect to $x \in \mathbb{R}^{n}$;
define accordingly the set $C^{k}(\mathbb{R}^{n};\mathbb{B})$.
%In case of real-valued functions, i.e. $\mathbb{B}=\mathbb{R}$, we shall simply write $C^{0,k}\big([0,T]\times\mathbb{R}^{n}\big)$ and $C^{k}(\mathbb{R}^{n})$, respectively.

For any $1\leq
p\leq\infty$, any
$\lambda\in\mathbb{R}$, any $0\leq t<T<\infty$, and any Banach space $\mathbb{X}$ with norm
$\|\cdot\|_{\mathbb{X}}$, %and any sub-$\sigma$-algebra $\mathcal{G}\subseteq\mathcal{F}$,
denote by
\begin{itemize}
%\item $\mathbb{L}^{p}_{\mathcal{G}}(\Omega,\mathbb{X})$ the
%set of all $\mathbb{X}$-valued, $\mathcal{G}$-measurable random
%variables $X$ such that
%$E\|X\|_{\mathbb{X}}^{p}<\infty\; ;$

%\item $\mathbb{L}^{p}_{\mathbb{F}}(t,T;\mathbb{X})$ the set
%of all $\mathbb{F}$-progressively measurable, $\mathbb{X}$-valued
%processes $X:[0,T]\times\Omega\rightarrow\mathbb{X}$ such
%that $\int_{t}^{T}\|X(s)\|_{\mathbb{X}}^{p}ds<\infty, \
%\text{a.s.} ;$

%\item $\mathbb{M}^{p}(t,T;\mathbb{X})\equiv\mathbb{L}^{p}_{\mathbb{F}}\big(t,T;\mathbb{L}^{p}(\Omega;\mathbb{X})\big)$ the set
%of all $\mathbb{F}$-progressively measurable, $\mathbb{X}$-valued
%processes $X:[0,T]\times\Omega\rightarrow\mathbb{X}$ such
%that $E\int_{t}^{T}\|X(s)\|_{\mathbb{X}}^{p}ds<\infty\; ;$

%\item $C_{\mathbb{F}}([t,T];\mathbb{X})$
%the set of all $\mathbb{F}$-adapted, $\mathbb{X}$-valued processes
%$X(\cdot,\omega):[0,T]\rightarrow\mathbb{X}$ that are continuous on $[t,T]$ for $P$-a.e.
%$\omega\in\Omega\;;$

\item $C_{\mathbb{F}}\big([t,T];\mathbb{L}^{p}(\Omega;\mathbb{X})\big)$ the set of all $\mathbb{F}$-adapted, $\mathbb{X}$-valued processes
$X(\cdot,\omega):[0,T]\rightarrow\mathbb{X}$ that are continuous on $[t,T]$ for $P$-a.e. $\omega\in\Omega$ and satisfy $E\|X(s)\|_{\mathbb{X}}^{p}<\infty \ \ \text{for \ every} \ s \in [t,T] \;;$
%\end{itemize}
%We define similarly the set
%$C_{\mathbb{F}}\big([t,T];\mathbb{L}^{p}(\Omega;\mathbb{X})\big)$. Furthermore, for any
%$\lambda\in\mathbb{R}$, %and $\mathbb{X}$ be a Banach space equipped with a norm $||\cdot||_{\mathbb{X}}$
%denote by
%\begin{itemize}
\item $\mathbb{M}_{\lambda}^p\big(t,T; \mathbb{X}\big)$ the set of all $\mathbb{F}$-adapted, $\mathbb{X}$-valued processes $X:[0,T]\times\Omega\rightarrow\mathbb{X}$ such that
 $||X||_{\lambda}^{p}\triangleq E\int_{t}^{T} e^{\lambda s} ||X(s)||_{\mathbb{X}}^p ds<\infty.$
\end{itemize}
Define similarly all the above sets for an infinite time-horizon, i.e., on the
non-negative real numbers $\mathbb{R}_{+}$.
% stand and define the space $\mathbb{M}_{\lambda}^p\big(0,\infty; \mathbb{X}\big)\equiv\mathbb{L}_{\lambda}^p\big(\mathbb{R}_{+};\mathbb{L}^{p}(\Omega;\mathbb{X})\big)$.
 Let also  $\mathbb{R}^{+}$ stand for the positive real numbers.

For every random field $\mathrm{F}:[0,T]\times\mathbb{R}^{n}\times\Omega\rightarrow\mathbb{R}$ denote by $A^{\mathrm{F}}:[0,T]\times\mathbb{R}^{n}\times\Omega\rightarrow\mathbb{R}$
 and $\Psi^{\mathrm{F}}:[0,T]\times\mathbb{R}^{n}\times\Omega\rightarrow\mathbb{R}^{1\times d}$ the random fields for which $\mathrm{F}$ obtains the semimartingale decomposition
\[ \mathrm{F}(t,x)=\mathrm{F}(0,x)+\int_{0}^{t}A^{\mathrm{F}}(s,x)ds+\int_{0}^{t}\Psi^{\mathrm{F}}(s,x)dW(s),
 \ \ \ \forall\ \, (t,x) \in [0,T]\times\mathbb{R}^{n}, \]
 \vskip8pt
 \noindent
almost surely. Assume a similar notation both for stochastic processes and for an infinite time-horizon as well. By default in the sequel, all (in)equalities between random quantities are to be understood $dP$-almost surely, $dP\otimes dt$-almost everywhere or
$dP\otimes dt\otimes dx$-almost everywhere, as suitable in the situation at hand.

We state also next a useful implication
of the generalized IKW formula (e.g. Kunita \cite{Kunita},
Section 3.3, pp 92-93). This will enable us to carry out computations, regarding the change-of-variable formula for random fields,
in a stochastically modulated dynamic framework.

\begin{proposition}\label{prop:IKW}
Suppose that the random field
$\, \mathrm{F}:[0,T]\times\mathbb{R}^{n}\times\Omega\rightarrow\mathbb{R}$
is of class $C^{0,2}([0,T]\times\mathbb{R}^{n};\mathbb{R})$ and $A^{\mathrm{F}}:[0,T]\times\mathbb{R}^{n}\times\Omega\rightarrow\mathbb{R}$, $\Psi^{\mathrm{F}}:[0,T]\times\mathbb{R}^{n}\times\Omega\rightarrow\mathbb{R}^{1\times d}$ are  $\mathbb{F}$-adapted random fields of class  $C^{0,1}([0,T]\times\mathbb{R}^{n};\mathbb{R})$ and $C^{0,2}([0,T]\times\mathbb{R}^{n};\mathbb{R}^{1\times d})$, respectively. Furthermore,
let $\, X=(X^{(1)}, \ldots, X^{(n)})^{\top}\,$ be a
vector of continuous semimartingales, %with decompositions
%\[ X^{(i)}(t)=  X^{(i)}(0)+\int_{0}^{t}\mu^{(i)}(s)ds+\int_{0}^{t}\sigma^{(i)}(s)dW(s); \ \ \ i=1,\ldots,n, \]
where $\ A^{X^{(i)}}(\cdot)$ is an
almost surely integrable process and
$\Psi^{X^{(i)}}(\cdot)$ is
an $\mathbb{F}$-progressively measurable, almost surely square
integrable process, $i=1,2,...,n$. Then
$\mathrm{F}\big(\cdot,X(\cdot)\big)$ is also a continuous
semimartingale, with decomposition
\begin{align}
\mathrm{F}\big{(}t,X(t)\big{)}= \mathrm{F}\big{(}0,X(0)\big{)}
& + \int_{0}^{t}A^{\mathrm{F}}\big{(}s,X(s)\big{)}ds + \int_{0}^{t}\Psi^{\mathrm{F}}\big{(}s,X(s)\big{)}\;dW(s)
 \nonumber\\[10pt]
 &+\int_{0}^{t} \left\langle A^{X}(s),\nabla_x\mathrm{F}\big{(}s,X(s)\big{)}\right\rangle ds\nonumber\\[10pt]
 &+\int_{0}^{t}\left(\nabla_x\mathrm{F}\right)^{\top}\big{(}s,X(s)\big{)}\;\Psi^{X}(s)\;dW(s)
\nonumber
\end{align}
%\noindent
\begin{align}
%\label{2b}
& \quad + \frac{1}{2}\int_{0}^{t}tr\left\{\Big(\Psi^{X}\big(\Psi^{X}\big)^\top\Big)(s)\;\nabla_{xx}\mathrm{F}\big{(}s,X(s)\big{)}\,
\right\}ds\nonumber \\[10pt]
& \quad +\int_{0}^{t}tr\left\{\big(\nabla_x\Psi^{\mathrm{F}}\big)^{\top}\big{(}s,X(s)\big{)}\,\Psi^{X}(s)\right\}ds\nonumber
\end{align}
%\vskip4pt
\noindent
 for every $0\leq t\leq T$.
 \end{proposition}

\vskip10pt

\subsection{Krugman's target zones exchange rate model}\label{KRUGMAN-SECTION}

Krugman \cite{Krugman}
developed a  simple economic model for the exchange rate behavior subject to a target zone regime. According to this model,
centering the target zone around zero and expressing all variables in natural logarithms,
the only exogenous source of noise  is introduced by the \emph{velocity shift term} $v(\cdot)$ which evolves according to the FSDE
\begin{equation*}%\label{eqn:motFSDEkrug}
\begin{split}
dv(t)&=\sigma dW(t), \quad t> 0,\\[8pt]
v(0)&=v;
\end{split}
\end{equation*}
i.e., it is a Brownian motion with constant intensity $\sigma\neq 0$
and does not anticipate future changes. Then the \emph{foreign exchange rate} $s(\cdot),$ which is assumed to depend on the expected rate of future depreciation, experiences velocity shocks
according to

\begin{equation*}
s(t)=m(t)+v(t)+\gamma\frac{E\big[ds(t)\big|\mathcal{F}(t)\big]}{dt}, \quad t\geq 0,
\end{equation*}
\vskip10pt
\noindent
where $m(\cdot)$ is the \emph{domestic money supply} and the last term
represents the \emph{expected rate of depreciation}.  This term involves the value of the variable $s(\tau)$ at future times $\tau > t$, a fact made clearer by stating $s(\cdot)$ in the integrated form

\begin{equation*}%\label{eqn:motBSDEkrug}
s(t)=\frac{1}{\gamma}E\left[\int_{t}^{\infty}e^{-\frac{1}{\gamma}(\theta-t)}\big(m(\theta)+v(\theta)\big)d\theta\bigg|\mathcal{F}(t)\right],\quad t\geq 0,
\end{equation*}
\vskip5pt
\noindent
which represents the foreign exchange rate as a present discounted value of future realizations of $m(\cdot)+v(\cdot)$. Obviously, future expectations for the stochastic process $m(\cdot)+v(\cdot)$ drive the exchange variability today.

 It is clear that this model has an FBSDEs flavor; in fact, it is equivalent (see Yannacopoulos \cite{Yannacopoulos}) to an infinite horizon (decoupled) FBSDEs system  of the form
\begin{equation}\label{KRUGMAN-BSDE}
\begin{aligned}
dv(t)&=\sigma dW(t), \quad t> 0, \\[8pt]
ds(t)&=\left[-\frac{1}{\gamma}(m(t)+v(t)) + \frac{1}{\gamma} s(t)\right]dt -z(t) dW(t), \quad t\geq 0,\\[8pt]
v(0)&=v,
\end{aligned}
\end{equation}
\vskip5pt
\noindent
where $s(\cdot)$ satisfies a final condition so that it decays to $0$ as $t \to \infty$ and $z(\cdot)$ is a stochastic process (we skip technicalities for the time being concerning its regularity) that plays an essential role in ensuring the right asymptotic behavior for $s(\cdot)$.
Furthermore,
the money supply $m(\cdot)$ can be considered as an  exogenous policy variable  of the model, which can be chosen so as to keep the exchange rate $s(\cdot)$ within the target zone, thus turning the model into the form of a controlled BSDE (cf. Yannacopoulos \cite{Yannacopoulos05}). In this model $v(\cdot)$ plays the role of the fundamental and $s(\cdot)$ the role of the asset.

It is quite natural to assume that the discount factor $\gamma$ and the intensity $\sigma$ of the shocks in \eqref{KRUGMAN-BSDE} are not constant but depend on time in a random fashion, i.e., are stochastic processes $\gamma(\cdot)$ and $\sigma(\cdot)$. If we therefore assume the discount factor to vary randomly in time (as the discount factor can well be assumed to depend on the state of the world) and the velocity shocks to be subject to random intensities then \eqref{KRUGMAN-BSDE} becomes a (decoupled) BSDE with random coefficients, which may model more accurately and more robustly the fluctuations of the exchange rate around equilibrium and thus lead to a more efficient money supply policy for the confinement of the exchange rate in the target zone.

\vskip15pt

\subsection{Dornbusch's model of exchange rate overshooting}\label{DORNBUSCH-SECTION}

A stochastic extension of the overshooting Dornbusch model of exchange rate determination (see e.g. Dornbusch \cite{Dornbusch})
was investigated by Neely et al \cite{Neely-Weller-Corbae}. All the  variables introduced in the model  except for interest rates were expressed in natural logarithms.
 %\cite{Neely-Weller-Corbae}
The first equation of the model
\[m(t)-p(t) =\nu y(t)-\xi i(t), \quad t\geq 0,\]
is the equilibrium condition of the domestic money market, where $m(\cdot)$ is the \emph{domestic money supply},
 $p(\cdot)$ is the \emph{domestic price level}, $y(\cdot)$ is the \emph{level of output} in the economy and
 $i(\cdot)$ is the nominal \emph{domestic interest rate}.
Denoting by $s(\cdot)$ the \emph{domestic price of foreign exchange}, the second equation

\begin{equation}\label{motdorny}
y(t) = -\vartheta\left(i(t) -\frac{E\big[dp(t)\big|\mathcal{F}(t)\big]} {dt}\right) +\eta\big(s(t)-p(t)\big), \quad t\geq 0,
\end{equation}
\vskip8pt
\noindent
is the condition for a goods market equilibrium, where the level of output
is negatively correlated ($\vartheta>0$) with the \emph{real interest rate} $i(t) -E[dp(t)|\mathcal{F}(t)]/dt,\ t\geq 0,$
and positively correlated ($\eta>0$) with the \emph{real exchange rate}  $s(\cdot)-p(\cdot)$ that corresponds to competitiveness.
The third equation

\[\frac{E\big[ds(t)\big|\mathcal{F}(t)\big]} {dt} = i(t) - i^{\ast}, \quad t\geq 0,\]
\vskip8pt
\noindent
is an uncovered interest parity condition, in which the
expected rate of depreciation of the domestic currency equals the nominal
interest differential that could also incorporate the possibility of jumps as a consequence of a realignment.
The last equation of the model
\begin{equation}\label{motdorp}
\begin{split}
dp(t) &= \phi \big( y(t)-\bar{y}\big)dt +\sigma dW(t), \quad t> 0,\\[12pt]
p(0)&=p,
\end{split}
\end{equation}
represents the less-than-instantaneous price adjustment in terms of the output, in excess of
its long-run full employment level $\bar{y}$, and a stochastic term that captures the external random shocks in the economy.

Setting for matters of simplicity $i^{\ast}=\bar{y}=0,$ one may solve the linear system of the first two equations of the model
for $i(\cdot)$ and $y(\cdot)$, and substituting them in the last two equations arrive at the system
\vskip5pt
\noindent
\begin{align*}%\label{motps}
dp(t)&=\left[\frac{\phi\vartheta}{D}m(t)-\frac{\phi(\vartheta+\xi\eta)}{D}p(t)+\frac{\phi\xi\eta}{D}s(t)\right]dt +\sigma dW(t),\quad t>0, \nonumber\\[12pt]
s(t)&=\frac{1-\phi\vartheta}{\nu\eta}m(t)-\frac{1-\phi\vartheta-\nu\eta}{\nu\eta}p(t)+\frac{D}{\nu\eta}\cdot\frac{E\big[ds(t)\big | \mathcal{F}(t)\big]}{dt}, \quad t\geq 0,\\[12pt]
p(0)&=p,\nonumber
\end{align*}
\vskip5pt
\noindent
where $D\triangleq \nu\vartheta+\xi-\phi\vartheta\xi$. The latter equation of this system can be reformulated as

\[s(t)=\frac{\nu\eta}{D}E\left[\int_{t}^{\infty}e^{-\frac{\nu\eta}{D}(\theta-t)}\left(\frac{1-\phi\vartheta}{\nu\eta}m(\theta)-\frac{1-\phi\vartheta-\nu\eta}{\nu\eta}p(\theta)\right)d\theta\bigg|\mathcal{F}(t)\right],\]
\vskip8pt
\noindent
which combined with the former equation imposes a more general future expectations model than the one derived in Krugman model. This model has been expressed in terms of an infinite horizon FBSDEs system in Yannacopoulos \cite{Yannacopoulos} as

\begin{equation}\label{DORNBUSCH-FBSDE}
\begin{aligned}
dp(t)&=\frac{1}{D}\Big[\phi\vartheta m(t)-\phi (\vartheta + \xi \eta) p(t)+ \phi \xi \eta s(t)\Big]dt + \sigma dW(t),\quad t>0, \\[10pt]
ds(t)&=\frac{1}{D}\Big[(\phi\vartheta-1)m(t)+(1-\nu \eta -\phi \vartheta) p(t) + \nu \eta s(t)\Big] dt -z(t) dW(t),\quad t\geq 0,\\[12pt]
p(0)&=p,
\end{aligned}
\end{equation}
where $s(\cdot)$ has to satisfy an appropriate asymptotic condition as $t \to \infty$ and $z(\cdot)$ is an unknown stochastic process which essentially guarantees the right asymptotic behavior of $s(\cdot)$. Here $p(\cdot)$ plays the role of the fundamental variable and $s(\cdot)$  the role of the asset.

All the coefficients used in this model were assumed to be constants or at most deterministic functions.  This is clearly an oversimplification as the parameters involved in the model are subject to random  variability and even if they are assumed to be constant, the econometric techniques required to calibrate the model yield estimators that are random processes assumed to converge to the true value of the parameter. It is therefore important to
extend this model so as to accommodate more sources of randomness in order to enhance its descriptive properties against the real economy. In particular,
the price level $p(\cdot)$ in (\ref{motdorp}) may both undergo real exogenous shocks according to a stochastic volatility rate $\sigma(\cdot)$ and adjust in terms of the excess output level
with a random speed $\phi(\cdot).$ Moreover, by assuming stochastic coefficients $\vartheta(\cdot)$ and $\eta(\cdot)$ in (\ref{motdorny}) then the model could produce realistic patterns of correlation
between the output $y(\cdot)$ and exchange rates or interest rates differentials. This turns \eqref{DORNBUSCH-FBSDE} quite naturally to an infinite horizon random FBSDEs system that is expected to provide a more realistic model for the overshooting patterns observed in exchange rate dynamics. The study of  this system and its solvability will provide more realistic conditions on the restrictions needed in the system so as to achieve equilibrium in the long run.
\vskip10pt
\subsection{Blanchard's model for the interaction of the real economy with the stock market}

Blanchard \cite{Blanchard} proposed a model (presented here as revisited in Miller and Weller \cite{Miller-Weller}) to study the interaction between the output of the real economy
with the stock market, in an attempt to model the possible stabilization or destabilization effects of expectations concerning the real economy on the stock market.
The \emph{level of demand} $Y(\cdot)$ is expressed in terms of the \emph{output} $x(\cdot)$,
the \emph{stock market} value $s(\cdot)$ and
the \emph{index of fiscal policy} $g(\cdot)$  by
\begin{equation}\label{eqn:motd}
Y(t)=a_1 x(t)+a_2 s(t)+g(t), \quad t\geq 0,
\end{equation}
where $a_{1},a_{2}$ are considered appropriate constants.
Furthermore, output adjusts to demand and Brownian stochastic shocks over time, as postulated by the rule
\begin{align}\label{eqn:motx}
\begin{split}
dx(t)&= a_3\big(Y(t)-x(t)\big)dt+\sigma dW(t),\quad t>0,\\[8pt]
x(0)&=x,
\end{split}
\end{align}
where $a_{3}$ is an appropriate proportionality factor.
Finally, in order to exclude arbitrage the rate of return on stock market shares
should be equal to the yield on short-term bonds, subject to the present of a risk premium
captured by the consumption capital asset price model; that is
\begin{equation}\label{eqn:motq}
\frac{E\big[ds(t)\big|\mathcal{F}(t)\big]} {dt} = \big(i(t)+\rho\, \vartheta_{s}\,\sigma^2\big)s(t)-c(t), \quad t\geq 0,
\end{equation}
where $\rho$ is an index of risk aversion, $\theta_{s}\sigma^2$ reflects the covariance between aggregate consumption and the stock market, and $i(\cdot)$ is the short
term \emph{nominal interest rate}. In Blanchard \cite{Blanchard} and Miller and Weller \cite{Miller-Weller} feedback stabilization rules on the consumption and the fiscal policy turn $c(\cdot)$ and $g(\cdot)$
to linear functions of the output $x(\cdot)$ and thus  close the system.

The equations of \eqref{eqn:motd}-\eqref{eqn:motq} can be expressed as an infinite horizon FBSDEs system of the form
\begin{equation}\label{BLANCHARD-FBSDE}
\begin{aligned}
dx(t)&= \Big[a_3 (a_1 -1) x(t)+ a_{2} a_{3} s(t)+ a_{3} g(t) \Big]dt+\sigma dW(t), \quad t>0, \\[10pt]
ds(t)&=\Big[( i(t) +\rho\, \vartheta_{s}\,\sigma^2\big)s(t)-c(t)\Big] dt - z(t)dW(t), \quad t\geq 0,\\[10pt]
x(0)&=x,
\end{aligned}
\end{equation}
where  $s(\cdot)$ satisfies a proper asymptotic condition as $t\to \infty$. The stochastic process $z(\cdot)$ is the auxiliary process needed
for the well-posedness of the system and $g(\cdot)$, $c(\cdot)$ can be considered as control variables chosen optimally so that  the system is driven to the desired state (rather than set a priori to a feedback rule as in the existing literature without any reference at to its optimality).

 The various parameters in system \eqref{BLANCHARD-FBSDE} are subject to random fluctuations, and especially  $i(\cdot)$. This naturally leads to considering  \eqref{BLANCHARD-FBSDE}
as a random coefficient controlled infinite horizon FBSDEs system which can better capture the adjustment of prices over time to equilibrium  if the model is well-posed or the instability caused by expectations otherwise. Furthermore, an optimal control approach  of the resulting random FBSDEs may serve in the efficient design of an economic policy rule that will be more robust and perform better than the ad hoc feedback rules often employed in practice.

\vskip10pt
\subsection{Black's consol rate conjecture}\label{sub:consolrate}

Another interesting example of future expectations model which is important in the field of mathematical finance can be found in Duffie et al \cite{Duffie-Ma-Yong} and is related to the famous Black's
consol rate conjecture connecting the short rate process with the price of long term bonds and in particular perpetuities.  Should the conjecture hold, the \emph{short rate} process $r(\cdot)$ and the \emph{price} process  $Y(\cdot)$ of perpetuities paying dividends at infinitum at the constant rate $1$  (\emph{consols}) should be related via the following pattern:
\begin{align*}
dr(t)&=\mu\big(r(t), Y(t)\big)dt+\alpha\big(r(t),Y(t)\big)dW_0(t),\quad t>0,\\[12pt]
Y(t)&=E_0\left[\int_{t}^{\infty}\exp\left(-\int_{t}^{s}r(u)du\right)ds\Big|\mathcal{F}(t)\right],\quad t\geq 0,\\[12pt]
r(0)&=r>0,
\end{align*}
where
$\mu:\mathbb{R}_{+}\times\mathbb{R}^{+}\rightarrow\mathbb{R}$ and $\alpha:\mathbb{R}_{+}\times\mathbb{R}^{+}\rightarrow\mathbb{R}^d$ are given functions that satisfy certain technical conditions,  and
$E_0$ denotes expectation under the equivalent martingale measure that makes the $d$-dimensional process $W_0(\cdot)$ a standard Brownian motion.

This model connects future expectations concerning the short rate with the bond market and clearly is an FBSDEs system, whose well-posedness may yield results on the equilibrium of the bond market. Duffie et al \cite{Duffie-Ma-Yong} studied this problem using FBSDEs and thus proved Black's conjecture under rather general conditions, as long as the relevant functions are deterministic. In fact, they also showed that a deterministic functional relation should hold between $Y(\cdot)$ and $r(\cdot)$ which guarantees Markovian dynamics, and furthermore, the representing function is characterized by the solution of a
deterministic elliptic PDE.

 Once again we could postulate that the diffusion functions
$\mu$ and $\alpha$ are now random fields, so as to generalize the influence the consol rate $Y(\cdot)$ has on the stochastic behavior of the short rate $r(\cdot)$.

\vskip10pt

\subsection{A general future expectations model and its connection with infinite horizon random FBSDEs}

Motivated by the examples presented in Subsections \ref{KRUGMAN-SECTION}-\ref{sub:consolrate} we formulate a general model for the interaction between
future expectations and present actions of the form:
\begin{align}\label{REM}
dX(t)&=b\big(t,X(t),Y(t)\big)dt+\sigma\big(t, X(t), Y(t)\big)dW(t), \quad t>0,\nonumber\\[13pt]
Y(t)&=E\left[\int_{t}^{\infty}e^{-\int_{t}^{s}r(\theta)d\theta}\;g\big(s,X(s), Y(s)\big)ds\bigg|\mathcal{F}(t)\right], \quad t\geq 0,\\[13pt]
X(0)&=x\nonumber
\end{align}
for $\mathbb{F}$-adapted stochastic processes $X:\mathbb{R}_{+}\times \Omega \rightarrow  \mathbb{R}^{n}$ and  $Y:\mathbb{R}_{+}\times \Omega \rightarrow  \mathbb{R}^{m},$ where
$b:\mathbb{R}_{+}\times \mathbb{R}^{n}\times \mathbb{R}^{m}\times \Omega\rightarrow  \mathbb{R}^{n},$ $\sigma :\mathbb{R}_{+}\times  \mathbb{R}^{n}\times \mathbb{R}^{m}\times \Omega\rightarrow  \mathbb{R}^{n\times d},$
$g:\mathbb{R}_{+}\times \mathbb{R}^{n}\times \mathbb{R}^{m}\times \Omega\rightarrow  \mathbb{R}^{m}$ and $r:\mathbb{R}_{+}\times \Omega \rightarrow  \mathbb{R}$ are given random fields (and not deterministic functions). In this general formulation of future expectations systems, $X(\cdot)$ plays the role of the \emph{fundamental variables} and  $Y(\cdot)$ the role of the\emph{ assets}. This is a general \emph{stochastic saddle point system}, whose linear form with constant coefficients  was proposed in Miller and Weller \cite{Miller-Weller} and whose nonlinear version but with deterministic coefficients was studied in detail in Yannacopoulos \cite{Yannacopoulos}. The introduction of random coefficients is  by no means a trivial generalization and is important from the point of view of economic modeling, therefore, the future expectations system \eqref{REM}
will be the main object of the present chapter.
Stochastic control aspects of this system %\eqref{REM}
 are investigated in detail in Section \ref{sec:ihmp}, where we incorporate in the coefficients of the system
a control process $ u(\cdot)$ that can be chosen so as to drive the latter to a desired state in the long run.

 In order to ensure that the variables of the above system grow at a reasonable rate, leading to a stable solution, when any of them diverge from the equilibrium states,
 we shall need to make certain  hypotheses about the growth and regularity of the aforementioned random fields; a rigorous definition for the solution of (\ref{REM}) is given in Definition \ref{DEF:SOLREM}.
 %Towards this direction, we shall need first some basic notation.

To prevent repetition, we shall now make assumptions concerning only the coefficients of the backward component of the future expectations model (\ref{REM}),
while its forward component will be treated in a more general framework in the next section; see Hypothesis \ref{Hyp} and Remark \ref{rem:fbsderem}.

\begin{hypothesis}\label{hyprem}
The coefficients of the BSDE of the saddlepoint system (\ref{REM}) satisfy the following conditions.

(H1) The stochastic discount factor $r:\mathbb{R}_{+}\times \Omega\rightarrow \mathbb{R}$ is $\mathbb{F}$-progressively measurable
 and uniformly bounded, i.e., there exists $\rho_{0}>0$ such that
 \[|r(\cdot)|\leq \rho_{0}.\]

(H2) The random field $g:\mathbb{R}_{+}\times \mathbb{R}^{n}\times \mathbb{R}^{m}\times \Omega \rightarrow \mathbb{R}^m $ is $\mathbb{F}$-progressively measurable for every $(x, y) \in \mathbb{R}^{n}\times \mathbb{R}^{m}$ and is continuous with respect to these variables. In addition, it is  monotonous with respect to $y,$ uniformly Lipschitz continuous with respect to $x$  and at most linearly growing in $y$; i.e., there exist $\mu_{0} \in \mathbb{R}$ and $c_0\geq 0$ such that
for every $t,\, x, \, y, \,  x_{i}, \, y_{i}, \ i=1,2,$ it satisfies
\begin{align*}
\langle g(t,x,y_1)-g(t,x,y_2), y_1-y_2\rangle &\leq \mu_0|y_{1}-y_{2}|^2,\\[8pt]
|g(t,x_1,y)-g(t,x_2,y)|&\leq c_0|x_{1}-x_{2}|, \\[8pt]
| g(t,0,y)|&\leq |g(t,0,0)|+\varphi\big(|y|\big)
\end{align*}
for a continuous increasing function  $\varphi: \mathbb{R}_{+}\rightarrow \mathbb{R}_{+},$ respectively. Furthermore, we have that
\[g(\cdot,0,0)\in \mathbb{M}_{\lambda}^2\big(0,\infty; \mathbb{R}^m\big).\]
\end{hypothesis}

\begin{definition} \label{DEF:SOLREM}
The adapted stochastic pair process  $(X,Y)$  constitutes a solution of the stochastic saddlepoint system (\ref{REM}) if it satisfies
its relationships and belongs to $\mathbb{M}_{\lambda}^2\big(0,\infty; \mathbb{R}^n\times\mathbb{R}^m\big)$
 for some $\lambda\in\mathbb{R}$.
  \end{definition}

The membership of $(X,Y)\in\mathbb{M}_{\lambda}^2\big(0,\infty; \mathbb{R}^n\times\mathbb{R}^m\big)$ is in fact a \emph{transversality condition} aiming to the well-posedness of the stochastic saddlepoint system  (\ref{REM}), subject to an appropriate adjusting rate $\lambda$.
It ensures that its adapted solution $\big(X(t), Y(t)\big),$ $t\geq 0,$  does not blow up as $t$ approaches infinity and eventually
is driven to equilibrium after long enough time.
The following result formulates the infinite horizon system of FBSDEs that is equivalent to the future expectations model (\ref{REM}).
This is an extension of Proposition 2.1 in Yannacopoulos \cite{Yannacopoulos}
%\cite{Yannacopoulos}
for generalized random coefficients and its proof, which we omit, follows through a similar methodology; cf. Appendix A.1 and Remark 2.1 therein.
\begin{proposition}  \label{Prop:FBSDEsToREM}
 Under Hypothesis \ref{hyprem}, the stochastic saddlepoint system (\ref{REM}), subject to adapted solutions in the sense of Definition \ref{DEF:SOLREM}, is equivalent to the infinite horizon FBSDEs system:
\begin{align}\label{FBSDEREM}
dX(t)&=b\big(t,X(t),Y(t)\big)dt+\sigma\big(t,X(t),Y(t)\big)dW(t), \quad t> 0,\nonumber\\[10pt]
dY(t)&=\left[-g\big(t, X(t),Y(t)\big)+r(t)Y(t)\right]dt+Z(t)dW(t), \quad t\geq 0,\\[10pt]
X(0)&=x\in\mathbb{R}^{n},\nonumber
\end{align}
 where
 $(X, Y, Z)\in \mathbb{M}_{\lambda}^2\big(0,\infty; \mathbb{R}^n\times\mathbb{R}^m\times \mathbb{R}^{m\times d}\big)$  for some $\mathbb{F}$-adapted stochastic process  $Z(\cdot)$ .
 %Furthermore, if (H3) and (H4) of Hypothesis \ref{Hyp} hold with $c_2=k_2=k_5=0,$ i.e.,
 %the random fields $b,$ $\sigma$ and $f$ are independent of the variable $z,$ then the converse holds as well
%and the two models are equivalent.
%the infinite horizon system of FBSDEs (\ref{GFBSDEs}), subject to adapted solutions in the sense of Definition \ref{DEF:SOLINFBSDE}, is equivalent to the rational expectations model (\ref{REM}).
\end{proposition}
\begin{remark}
As a direct consequence of the preceding proposition, the well-posedness of the stochastic saddlepoint system (\ref{REM}) reduces to the solvability of the equivalent infinite horizon random FBSDEs system (\ref{FBSDEREM}).
This is studied in the sequel and in the proper functional setting allows us to obtain relevant information concerning the effect of future expectations to the
evolution of the system, whether they are expected to stabilize or destabilize it and to what extent. The economic interpretation of $Z(\cdot)$ is the price of the
control that has to be enforced to this system in order to eventually reach the equilibrium state.
Technically, it is the auxiliary process that makes the backward component $Y(\cdot)$ of the FBSDEs (\ref{FBSDEREM}) to be adapted at all times. %The next result
\end{remark}
%\begin{proposition}\label{Prop:FBSDEsEquivRem}
%Let (H3) and (H4) of Hypothesis \ref{Hyp} hold with $c_2=k_2=k_5=0$; that is, the random fields $b,$ $\sigma$ and $f$ are independent of the variable $z.$
%Then the infinite horizon system of FBSDEs (\ref{GFBSDEs}), subject to adapted solutions in the sense of Definition \ref{DEF:SOLINFBSDE}, is equivalent to the rational expectations model (\ref{REM}).
%\end{proposition}

% Furthermore,
 %by extending the four-step-scheme using backward stochastic parabolic equations (rather than deterministic elliptic equations)
  %we may provide in Section \ref{sub:svs} a functional relationship between the variables $X$ and $Y$ in terms of a random functional, which will lead to a new notion of stochastic viscosity solutions.
   %Finally by formulating and studying properly selected optimal control problems for the random FBSDEs in Section \ref{sec:ihmp}, we obtain optimal intervention policies for the system which may drive %
   %it to a desired equilibrium state.
\begin{remark}
Regarding the solvability of the infinite horizon FBSDEs system (\ref{FBSDEREM}),
%it would be natural to try and employ the results of the previous chapter where, through the four-step-scheme for random variables
 %(cf. Ma and Yong \cite{Ma-Yong98}), in the one-dimensional case this system may obtain a solution in terms of a linearizable Burgers BSPDE; see Section \ref{sec:BRBE}.
 we first consider
 %shall explore this direction by considering
 for $n=m=d=1$ the infinite horizon \emph{linear} random FBSDEs system
\begin{align*}
  dX(t)&=\Big[\alpha(t) X(t)+\beta(t)Y(t)\Big]dt+\sigma(t)dW(t),\quad t>0,\nonumber\\[10pt]
  dY(t)&=\Big[\gamma(t)X(t)+\delta(t)Y(t)\Big]dt+Z(t)dW(t),\quad t\geq0,\nonumber\\[10pt]
  X(0)&=x\in\mathbb{R}, %\qquad Y(T)=p\big(X(T)\big),
  \end{align*}
where $\alpha, \beta, \gamma, \delta: \mathbb{R}_{+}\times \Omega\rightarrow \mathbb{R}$ are bounded stochastic processes% and $p:\Omega\rightarrow \mathbb{R}$ is a bounded $\mathcal{F}(T)$-adapted random variable
, which is the simplest possible form the system of
(\ref{FBSDEREM})
may take so as its random coefficients satisfy the conditions of Hypothesis \ref{hyprem}. Given a finite time-horizon and a linear terminal condition, Yong \cite{Yong06} introduced a decoupling
reduction method in order to obtain the solvability for a very big class of such linear FBSDEs.

 Alternatively, inspired by the four-step-scheme of Ma and Yong \cite{Ma-Yong98},
if one is looking for solutions of the form $Y(t)=U\big(t,X(t)\big)$, $t\geq0$,
the IKW formula yields that the pair of random fields  $(U,\Psi^{U})$ % is of class $ C_{\mathbb{F}}\big([0,T];\mathbb{L}^2(\Omega;C^3(\mathbb{R}))\big)\times \mathbb{L}_{\mathbb{F}}^2\big(0,T;C^2(\mathbb{R})\big)$
should be a solution of
% solves
the following  Burgers type infinite horizon BSPDE
\begin{align*}
 dU&=\Big[-\frac{1}{2}\;    \sigma^2(t)U_{xx}-\beta(t)UU_{x}-\alpha(t)xU_x-\sigma(t)\Psi_{x}^{U}+\delta(t)U+\gamma(t)x\Big]dt\\[8pt]
 &\quad+\Psi^{U}dW(t), \quad t\geq0. \nonumber
%U(T,x)&=p(x),\quad x\in \mathbb{R}.
 \end{align*}
In case of a linear random field $U$ for a finite time-horizon, the solvability of this BSPDE reduces to the solvability of a stochastic Riccati BSDE that was studied by Yong \cite{Yong06}.
However if $U$ is not linear, apart from the infinite time-horizon, this type of Burgers equation is not of the same form with the one studied by Englezos et al. \cite{Englezosetal}, and accordingly we can not make use of
  a generalized version of the Cole-Hopf transformation in order to linearize it and obtain solutions in terms of a backward stochastic heat equation.

   Consequently, the section that follows deals with the solvability of the infinite horizon FBSDEs system (\ref{FBSDEREM}), whose random coefficients are not necessarily linear. A further investigation on the
   four-step-scheme and on the solvability of a resulting more general version of the above BSPDE takes place in Section \ref{sub:svs}.
\end{remark}

\vskip20pt

\section{Solvability And Sensitivity Analysis Of Infinite Horizon Random FBSDEs}\label{sec:ssa}

In this section, by means of a fixed-point argument, we shall establish existence and uniqueness for an adapted solution, i.e. well-posedness, of a more general infinite horizon FBSDEs system than the one of (\ref{FBSDEREM})
that has the form

\begin{align}\label{GFBSDEs}
dX(t)&=b\big(t,X(t),Y(t),Z(t)\big)dt+\sigma\big(t,X(t),Y(t),Z(t)\big)dW(t), \quad t> 0,\nonumber\\[10pt]
dY(t)&=-f\big(t,X(t),Y(t),Z(t)\big)dt+Z(t)dW(t), \quad t\geq0,\\[10pt]
X(0)&=x\in\mathbb{R}^n ,\nonumber
\end{align}
where $b, \, \sigma$ and $f$ are random fields with properties determined by the subsequent assumptions.
\vskip10pt
\subsection{ Definitions and assumptions}

We make first the following assumptions.

\begin{hypothesis}\label{Hyp}
The coefficients of FBSDEs (\ref{GFBSDEs}) satisfy the following conditions.

(H1) The random fields $b: \mathbb{R}_{+}\times \mathbb{R}^{n}\times \mathbb{R}^{m}\times \mathbb{R}^{m\times d}\times \Omega \rightarrow \mathbb{R}^n , \ \sigma:\mathbb{R}_{+}\times \mathbb{R}^{n}\times \mathbb{R}^{m}\times \mathbb{R}^{m\times d}\times \Omega \rightarrow \mathbb{R}^{n\times d}$ and  $f:\mathbb{R}_{+}\times \mathbb{R}^{n}\times \mathbb{R}^{m}\times \mathbb{R}^{m\times d}\times \Omega \rightarrow \mathbb{R}^m $ are $\mathbb{F}$-progressively measurable for every $(x, y, z) \in \mathbb{R}^{n}\times \mathbb{R}^{m}\times \mathbb{R}^{m\times d}$. In addition, they are continuous with respect to these variables.

(H2) There exist $\mu_{1}, \mu_{2}\in \mathbb{R}$ such that for every $t, \, x,  \,  y,  \, z,  \, x_{i}, \, y_{i}, \ i=1,2,$ the random fields $b$ and $f$ satisfy  the monotonicity conditions:
\begin{align*}
\langle b(t,x_{1},y,z)-b(t,x_{2},y,z), x_{1}-x_{2}\rangle&\leq \mu_{1}|x_{1}-x_{2}|^2,\\[7pt]
\langle f(t,x,y_1,z)-f(t,x,y_2,z), y_{1}-y_{2}\rangle&\leq \mu_{2}|y_{1}-y_{2}|^2.
\end{align*}

 (H3) The random field $b$ is uniformly Lipschitz continuous with respect to $(y,z),$ at most linearly growing in $x,$ and $f$ is uniformly Lipschitz continuous with respect to $(x,z),$ at most increasingly growing in $y$.
 In particular, there exist $k, k_{i}, c_{i}\geq 0, \ i=1,2,$ such that for every $t, \, x,  \, y,  \, x_{i},  \, y_{i},  \, z _{i}, \ i=1,2, $ we have  that
\begin{align*}
| b(t,x,y_1,z_1)-b(t,x,y_2,z_2)|&\leq k_{1}|y_{1}-y_{2}|+k_2||z_{1}-z_{2}||,\\[7pt]
| f(t,x_1,y,z_1)-f(t,x_2,y,z_2)|&\leq c_{1}|x_{1}-x_{2}|+c_2||z_{1}-z_{2}||,\\[7pt]
| b(t,x,0,0)|&\leq |b(t,0,0,0)|+k\big(1+|x|\big),\\[7pt]
| f(t,0,y,0)|&\leq |f(t,0,0,0)|+\psi\big(|y|\big)
\end{align*}
for a continuous increasing function  $\psi: \mathbb{R}_{+}\rightarrow \mathbb{R}_{+}.$

(H4) The random field $\sigma$ is uniformly Lipschitz continuous with respect to $(x, y, z).$ This means that there exist $k_i\geq 0, \ i=3,4,5,$ such that for every $t,  \, x_i,  \, y_i,  \, z_i, \ i=1,2,$ we have  that
\begin{equation*}
||\sigma(t,x_1,y_1,z_1)-\sigma(t,x_2,y_2,z_2)||^2\leq k_3|x_1-x_2|^2+k_4|y_1-y_2|^2+k_5||z_1-z_2||^2.
\end{equation*}

(H5) There exist constants $\lambda\in \mathbb{R},$ $\epsilon_{i}>0, \ i=1,2 \, $ and $ \, \rho_{i}>0, \ i=1,2,\,$ such that
\begin{equation}\label{INEQcrho}
1-c_2\rho_{2}^{-1}-\gamma c_1\rho_{1}^{-1}>0, \qquad  c_2\rho_2^{-1}<1,
\end{equation}
and
\begin{align}\label{INEQL}
2\mu_2+c_1\rho_1+c_2\rho_2+\frac{c_1\rho_{1}^{-1}(1-c_2\rho_{2}^{-1})}{1-c_2\rho_{2}^{-1}-\gamma c_1\rho_{1}^{-1} }< \lambda < &-2\mu_1-k_3-k_1\epsilon_{1}-k_2\epsilon_{2}\\[10pt]
&-\big(k_1\epsilon_{1}^{-1}+k_4\big)\vee \big(k_2\epsilon_{2}^{-1}+k_5\big),\nonumber
\end{align}
\vskip5pt
\noindent
where $$0\leq\gamma\triangleq\frac{k_2\epsilon_{2}^{-1}+k_5}{\big(k_1\epsilon_{1}^{-1}+k_4\big)\vee \big(k_2\epsilon_{2}^{-1}+k_5\big)}\leq 1.$$
\vskip10pt
\noindent
In addition the following holds:
\[ \big(b(\cdot,0,0,0),f(\cdot,0,0,0),\sigma(\cdot,0,0,0)\big)\in \mathbb{M}_{\lambda}^2\big(0,\infty; \mathbb{R}^n\times\mathbb{R}^m\times \mathbb{R}^{n\times d}\big). \]
%\begin{align*}
%E\left[\int_{0}^{\infty} e^{\lambda t} \big(|b(t,0,0,0)|^2+||\sigma(t,0,0,0)||^2\big)dt\right]&<\infty\\
%E\left[\int_{0}^{\infty} e^{\lambda t} |f(t,0,0,0)|^2 dt\right]&<\infty.
%\end{align*}
\end{hypothesis}

The reader should observe that $\gamma$ in the preceding hypothesis is not well-defined only in the special case of
$k_1=k_2=k_4=k_5=0.$ However, this case corresponds to a decoupled infinite horizon system of FBSDEs (\ref{GFBSDEs}) since the forward equation does not
depend on the components of the backward one. Therefore, this system can be treated with simpler methods than the ones we shall present in this section. The remark that follows
describes explicitly the connection between the two infinite horizon FBSDEs systems of (\ref{FBSDEREM}) and (\ref{GFBSDEs}), subject to their associated Hypotheses \ref{hyprem} and \ref{Hyp},
respectively.

\begin{remark}\label{rem:lowerboundlambda}
Treating the lower bound of $\lambda$ in ($\ref{INEQL}$) as a two-variable function with respect to $(\rho_1,\rho_2)\in\mathbb{R}^{+}\times \mathbb{R}^{+}$, it turns out through standard multi-variable analysis that it obtains a global minimum value
%of $2\mu_2+2c_1+2c_1c_2\sqrt{\gamma}+c_1^2\gamma+ c_2^2$
for the critical points $\rho_1^{\ast}=c_1\gamma+c_2\sqrt{\gamma}+1$ and $\rho_2^{\ast}=c_2+c_1\sqrt{\gamma}$, which satisfy also the constraints of (\ref{INEQcrho}). Then, given any $\epsilon_i>0, \,i=1,2,$ the widest possible range for $\lambda$ implied by ($\ref{INEQL}$), as far as its lower bound is concerned, is now given by the relationship:
\begin{align}\label{WIDER-INEQL}
2\mu_2+2c_1+2c_1c_2\sqrt{\gamma}+c_1^2\gamma+ c_2^2< \lambda < &-2\mu_1-k_3-k_1\epsilon_{1}-k_2\epsilon_{2}\\[10pt]
&-\big(k_1\epsilon_{1}^{-1}+k_4\big)\vee \big(k_2\epsilon_{2}^{-1}+k_5\big).\nonumber
\end{align}
\end{remark}
\begin{remark}\label{rem:fbsderem}
Firstly, it is clear that if the random fields $b, \, \sigma$ and $f$ of (\ref{GFBSDEs}) are independent of the argument $z$ and $f(t,x,y)=g(t,x,y)-r(t)y$ then we retrieve the infinite horizon FBSDEs system
(\ref{FBSDEREM}). Secondly, the random fields $b$ and $\sigma$ of the latter system comply with the conditions of Hypothesis \ref{Hyp}, which refers to the generalized system (\ref{GFBSDEs}), by taking $k_2=k_5=0;$ thus, from now onwards we shall accordingly impose these conditions on the random fields $b$ and $\sigma$ of both FBSDEs systems. Finally, due to the above special form of the random field $f,$ Hypothesis \ref{hyprem} is consistent with Hypothesis \ref{Hyp}; in particular, if the former holds then so does the latter, as far as the random field $f$ is concerned, for $\mu_2=\rho_0+\mu_0,$ $c_1=c_0,$  $c_2=0,$ and $\psi\big(|y|\big)=\varphi\big(|y|\big)+\rho_0|y|.$
\end{remark}
We may also give the definition for an adapted solution of the infinite horizon FBSDEs system (\ref{GFBSDEs}).
 \begin{definition} \label{DEF:SOLINFBSDE}
 A triplet of $\mathbb{F}$-adapted processes $(X, Y, Z)$ is an adapted solution of FBSDEs  (\ref{GFBSDEs}) if for every $t\geq0$ satisfies
 \begin{align*}
 X(t)&=x+\int_{0}^{t} b\big(s,X(s),Y(s),Z(s)\big)ds+\int_{0}^{t} \sigma\big(s,X(s),Y(s),Z(s)\big)dW(s),\\[12pt]
 Y(t)&=\int_{t}^{\infty} f\big(s,X(s),Y(s),Z(s)\big)ds-\int_{t}^{\infty} Z(s) dW(s)
 \end{align*}
 \vskip2pt
 \noindent
 and $(X, Y, Z)\in \mathbb{M}_{\lambda}^2\big(0,\infty; \mathbb{R}^n\times\mathbb{R}^m\times \mathbb{R}^{m\times d}\big),$
 where $\lambda$ satisfies (\ref{INEQL}).
 %Furthermore, we can easy presume that if the triplet $(X, Y, Z)$ is a solution of FBSDEs  (\ref{GFBSDEs}) then from Burkholder-Davis-Gundy inequality
%and the assumptions (H1)-(H5) we arrive at
%\begin{equation*}
%E\left[sup_{t\geq 0} e^{\lambda t}|X(t)|^2+sup_{t\geq 0} e^{\lambda t}|Y(t)|^2\right]<\infty.
%\end{equation*}
  \end{definition}

 %The above system is connected directly with the rational expectation models which are a special case of (\ref{GFBSDEs}) via a specific choice of the random fields  $b, \sigma$ and $f$ (see Section 2).
 In particular,  one may use the integrability of $Z(\cdot)$, the conditions (H1), (H3) and (H5) of Hypothesis \ref{Hyp}, and Burkholder-Davis-Gundy's inequality to obtain that
 %\vskip2pt
\begin{equation}\label{eqn:sup}
E\left(\sup_{t\geq 0}e^{\lambda\,t}|X(t)|^2+\sup_{t\geq 0}e^{\lambda\,t}|Y(t)|^2\right)<\infty.
\end{equation}

\vskip14pt
\subsection{Solvability of stochastic saddlepoint systems}\label{sub:sss}
In this subsection, we extend Yin's \cite{Yin} result on existence and uniqueness of solutions for the FBSDEs (\ref{GFBSDEs}). % by suggesting a more general lower bound for the parameter $\lambda$ in relationship (\ref{INEQL}).
In particular, Theorem \ref{thm:exuniqgeneral} below extends Theorem 3.1 in Yin
\cite{Yin}
by providing the same upper bound but a more general lower bound for $\lambda$ via (\ref{INEQL}), subject to which the infinite horizon FBSDEs system of (\ref{GFBSDEs}) admits a unique and adapted solution. To see this,
according to Remark \ref{rem:lowerboundlambda} the best possible lower bound allowed for $\lambda$ is given by (\ref{WIDER-INEQL}), which thanks to $0\leq \gamma\leq 1$ is strictly smaller than $2\mu_2+4c_1^2+2c_2^2+1$
 that was proposed by Yin \cite{Yin}.% whenever we have that $c_1>0$ and/or $c_2>0$ then it suffices to take $\rho_1\triangleq4c_1>0$ and/or $\rho_2\triangleq2c_2>0$, otherwise the choice of
% $\rho_1$ and/or $\rho_2$ does not play any role, respectively, and recall that $0\leq \gamma\leq 1$ in order to verify that
%the inequalities of (\ref{INEQcrho}) hold and the lower bound of $\lambda$ in (\ref{INEQL}) becomes strictly smaller than $2\mu_2+4c_1^2+2c_2^2+1$, which was the one proposed by Yin \cite{Yin}. %Remark \ref{rem:yin},
%and proceed further by studying a comparison property of adapted solutions and their continuous dependence on a parameter.

First, we quote two lemmata which will be used in the subsequent analysis. The first one constitutes a direct consequence of Theorem 4.1 in Pardoux \cite{Pardoux}
and the other one is Lemma 3.2 in Yin \cite{Yin}.
\begin{lemma}\label{lem:existBSDE}
(Pardoux \cite{Pardoux}) Under %(H1), (H2), (H3) and (H5) of
Hypothesis \ref{Hyp} as regards the random field $f$ and given a process $\bar{X}(\cdot)\in\mathbb{M}_{\lambda_1}^{2}(0,\infty;\mathbb{R}^n)$ for $\lambda_{1}>2\mu_2+c_{2}^2,$
the BSDE
\begin{equation}\label{BSDE}
dY(t)=- f\big(t,\bar{X}(t),Y(t),Z(t)\big) dt-Z(t)dW(t), \qquad t\geq 0,
\end{equation}
admits a unique solution pair $(Y,Z)\in\mathbb{M}_{\lambda_{1}}^{2}(0,\infty;\mathbb{R}^m\times\mathbb{R}^{m\times d}).$
\end{lemma}

\begin{lemma}\label{lem:existFSDE}
(Yin \cite{Yin}) Under Hypothesis \ref{Hyp} as regards the random fields $b$, $\sigma$ and given a process pair $(\bar{Y},\bar{Z})$ that belongs to the class of
 $\mathbb{M}_{\lambda_{2}}^{2}(0,\infty;\mathbb{R}^m\times\mathbb{R}^{m\times d})$ for $\lambda_{2}<-2\mu_1-k_{3},$  the FSDE
\begin{align}\label{FSDE}
\begin{split}
dX(t)&= b\big(t,X(t),\bar{Y}(t),\bar{Z}(t)\big) dt-\sigma\big(t,X(t),\bar{Y}(t),\bar{Z}(t)\big)dW(t), \quad t> 0,\\[8pt]
X(0)&=x\in\mathbb{R}^n,
\end{split}
\end{align}
obtains a unique solution $X(\cdot)\in\mathbb{M}_{\lambda_{2}}^{2}(0,\infty;\mathbb{R}^n)$.
\end{lemma}

Clearly if $\lambda$ satisfies (\ref{INEQL}) then it also satisfies the inequalities of the above lemmata. Let Hypothesis \ref{Hyp} hold as well. Then for any given process
 $\bar{X}(\cdot)\in\mathbb{M}_{\lambda}^{2}(0,\infty;\mathbb{R}^n)$, due to Lemma \ref{lem:existBSDE}, there exists a unique process pair  $(Y,Z)\in\mathbb{M}_{\lambda}^{2}(0,\infty;\mathbb{R}^m\times \mathbb{R}^{m\times d})$ that solves the BSDE (\ref{BSDE}), and for this pair Lemma \ref{lem:existFSDE} guarantees a unique solution $X(\cdot)\in\mathbb{M}_{\lambda}^{2}(0,\infty;\mathbb{R}^n)$
 for the FSDE (\ref{FSDE}); hence, we may define accordingly the endomorphism $\Gamma_1: \mathbb{M}_{\lambda}^2(0,\infty;\mathbb{R}^n)\rightarrow \mathbb{M}_{\lambda}^2(0,\infty;\mathbb{R}^n):\bar{X}(\cdot)\mapsto X(\cdot).$  By the same token, we may define $\Gamma_2: \mathbb{M}_{\lambda}^{2}(0,\infty;\mathbb{R}^m\times\mathbb{R}^{m\times d})\rightarrow \mathbb{M}_{\lambda}^{2}(0,\infty;\mathbb{R}^m\times\mathbb{R}^{m\times d}):(\bar{Y},\bar{Z})\mapsto (Y, Z)$ to be  the endomorphism  which maps any given process pair $(\bar{Y},\bar{Z})$ to the unique solution $(Y,Z)$ of the BSDE (\ref{BSDE}), subject to  the unique solution $X(\cdot)$ of the FSDE (\ref{FSDE}).

 Given the above endomorphisms, we may now tackle the existence and uniqueness of an adapted solution for the strongly coupled system of FBSDEs (\ref{GFBSDEs}) by employing a \emph{fixed point scheme}.
 As we shall see in what follows, the inequality (\ref{INEQL})
 is a necessity for the parameter $\lambda$, since its lower and upper bounds not only ensure the existence of the backward and the forward component of FBSDEs (\ref{GFBSDEs})
 thanks to Lemmata \ref{lem:existBSDE} and \ref{lem:existFSDE}, respectively, but also enables the usage of the contraction mapping principle.

\begin{theorem}\label{thm:exuniqgeneral} Under Hypothesis \ref{Hyp} there exists a unique adapted solution $(X, Y, Z)$ to the FBSDEs (\ref{GFBSDEs}), in the sense of Definition \ref{DEF:SOLINFBSDE},  which belongs to the class of $\mathbb{M}_{\lambda}^2(0,\infty;\mathbb{R}^n\times\mathbb{R}^m\times \mathbb{R}^{m\times d}).$
\end{theorem}
\noindent
\textbf{Proof.} It is sufficient to show that the mapping $\Gamma_2$ is a contraction with respect to norm $||\cdot||_{\lambda}$, where $\lambda$ satisfies
(\ref{INEQL}). To this end, let $(\bar{Y}_i,\bar{Z}_i), \, i=1,2,$ be two process pairs in $\mathbb{M}_{\lambda}^2(0,\infty;\mathbb{R}^m\times \mathbb{R}^{m\times d})$  and from Lemma \ref{lem:existFSDE} consider the unique solutions
$X_i(\cdot), \, i=1,2,$ in $\mathbb{M}_{\lambda}^2(0,\infty;\mathbb{R}^n)$ of the corresponding FSDE (\ref{FSDE}), respectively.

A direct application of  It\^{o}'s rule, in conjunction with an elementary algebraic inequality and the conditions (H1)-(H4) of Hypothesis \ref{Hyp} for the random fields $b$ and $\sigma$, yields the relationship
 \begin{align*}\label{ineqx}
\bar{\lambda}_1 E\int_{0}^{T} e^{\lambda t}|X_1(t)-X_2(t)|^2dt&\leq  (k_1\epsilon_{1}^{-1}+k_4)E\int_{0}^{T} e^{\lambda t}|\bar{Y}_1(t)-\bar{Y}_2(t)|^2dt\nonumber\\ \\
&\quad \ +(k_2\epsilon_2^{-1}+k_5)E\int_{0}^{T} e^{\lambda t}||\bar{Z}_1(t)-\bar{Z}_2(t)||^2dt\nonumber
\end{align*}
for $T> 0$, where
 \begin{equation}\label{barl1}
 \bar{\lambda}_1\triangleq -\lambda-2\mu_1-k_3-k_1\epsilon_1-k_2\epsilon_2>0
\end{equation}
because of (\ref{INEQL}). Taking the limit $T\rightarrow \infty$ and employing Fatou's lemma we deduce that
\begin{equation}\label{metrx}
||X_{1}-X_{2}||_{\lambda}^2\leq \frac{k_1\epsilon_{1}^{-1}+k_4}{\bar{\lambda}_1}||\bar{Y}_1-\bar{Y}_2||_{\lambda}^2+\frac{(k_2\epsilon_2^{-1}+k_5)}{\bar{\lambda}_1}||\bar{Z}_1-\bar{Z}_2||_{\lambda}^2.
\end{equation}
\vskip5pt
Furthermore, according to Lemma \ref{lem:existBSDE} the corresponding  BSDE (\ref{BSDE}) given the above processes $X_i(\cdot), \, i=1,2,$ admits unique solution pairs
  $(Y_i,Z_i)$ of class $\mathbb{M}_{\lambda}^2(0,\infty;\mathbb{R}^m\times\mathbb{R}^{m\times d}), \, i=1,2,$ respectively. Then employing It\^{o}'s formula, an elementary algebraic inequality,
   the conditions (H1)-(H3) of Hypothesis \ref{Hyp} for the random field $f$, the inequalities (\ref{INEQcrho}) and (\ref{INEQL}), and Fatou's lemma we finally obtain
the inequality
\begin{align*}\label{ineqyz}
& e^{\lambda t}E| Y_{1}(t)-Y_{2}(t)|^2+\bar{\lambda}_2 E\int_{t}^{\infty} e^{\lambda s}|Y_{1}(s)-Y_{2}(s)|^2ds\nonumber\\[12pt]
&+(1-c_2\rho_{2}^{-1})E\int_{t}^{\infty} e^{\lambda s}||Z_{1}(s)-Z_{2}(s)||^2ds \leq \, c_1\rho_{1}^{-1}E\int_{t}^{\infty} e^{\lambda s}|X_{1}(s)-X_{2}(s)|^2ds\nonumber
\end{align*}
\vskip5pt
\noindent
for $t\geq 0$, where
 \begin{equation}\label{barl2}
 \bar{\lambda}_2\triangleq \lambda-2\mu_2-c_1\rho_1-c_2\rho_2>0.
\end{equation}
\vskip5pt
\noindent
This immediately yields the inequalities
\begin{align}\label{metryz}
||Y_{1}-Y_{2}||_{\lambda}^2&\leq \frac{c_1\rho_{1}^{-1}}{\bar{\lambda}_2}||X_{1}-X_{2}||_{\lambda}^2,\nonumber\\\\
||Z_{1}-Z_{2}||_{\lambda}^2&\leq \frac{c_1\rho_{1}^{-1}}{1-c_2\rho_2^{-1}}||X_{1}-X_{2}||_{\lambda}^2.\nonumber
\end{align}

Setting now
\begin{align}\label{defabc}
a&\triangleq ||X_1-X_2||_{\lambda}^2,\ \ \ \ \, \bar{a}\triangleq ||\bar{X}_1-\bar{X}_2||_{\lambda}^2,\nonumber\\[8pt]
b&\triangleq ||Y_1-Y_2||_{\lambda}^2,\quad \,\,\bar{b}\triangleq ||\bar{Y}_1-\bar{Y}_2||_{\lambda}^2, \\[8pt]
c&\triangleq ||Z_1-Z_2||_{\lambda}^2,\quad \,\,\bar{c}\triangleq ||\bar{Z}_1-\bar{Z}_2||_{\lambda}^2,\nonumber
\end{align}
the inequalities of (\ref{metrx}) and (\ref{metryz}) can be written as
\vskip5pt
\noindent
\begin{align}\label{ineqbarabc}
  a\leq \frac{1}{\bar{\lambda}_1}\Big[(k_1\epsilon_{1}^{-1}+k_4)\bar{b}+&(k_2\epsilon_2^{-1}+k_5)\bar{c}\Big], \quad b\leq\frac{c_1\rho_1^{-1}}{\bar{\lambda}_2}\,a,
\quad c\leq\frac{c_1\rho_1^{-1}} {1-c_2\rho_2^{-1}}\,a.
\end{align}
\vskip5pt
\noindent
Combine these inequalities and recall the parameter $\gamma$ in (H5) of Hypothesis \ref{Hyp} to get that
\begin{align*}
b+\gamma c\leq (\bar{b}+\gamma\bar{c})c_1\rho_1^{-1}\frac{\big(k_1\epsilon_{1}^{-1}+k_4\big)\vee \big(k_2\epsilon_{2}^{-1}+k_5\big)}{\bar{\lambda}_{1}}\left[\frac{1}{\bar{\lambda}_{2}}+\frac{\gamma}{(1-c_2\rho_2^{-1})}\right].
\end{align*}
\vskip5pt
\noindent
Under the relationships  (\ref{INEQcrho}), (\ref{INEQL}), (\ref{barl1}) and (\ref{barl2}) it is straightforward to see that
 \begin{align*}
 c_1\rho_1^{-1}\left[\frac{1}{\bar{\lambda}_{2}}+\frac{\gamma}{(1-c_2\rho_2^{-1})}\right]< 1 \qquad \text{and}\qquad
 \frac{\big(k_1\epsilon_{1}^{-1}+k_4\big)\vee \big(k_2\epsilon_{2}^{-1}+k_5\big)}{\bar{\lambda}_{1}}<1,
 \end{align*}
 \vskip5pt
\noindent
  which imply that the mapping $\Gamma_2 $ is a contraction subject to the norm $||\cdot||_{\lambda}$ and thus it has a unique fixed point.\hfill$\square$
\smallskip

\begin{remark} \label{rem:exisuniqsigma} %\label{thm:exisuniqsigma}
Note that Theorem \ref{thm:exuniqgeneral} holds in the special case of $k_5=0$ in Hypothesis \ref{Hyp}; that is, the diffusion coefficient $\sigma$ is independent of the argument $z$ and the system of FBSDEs (\ref{GFBSDEs})
reduces to
\begin{align}\label{FBSDEssigma}
dX(t)&=b\big(t,X(t),Y(t),Z(t)\big)dt+\sigma\big(t,X(t),Y(t)\big)dW(t), \quad t> 0,\nonumber\\[8pt]
dY(t)&=-f\big(t,X(t),Y(t),Z(t)\big)dt+Z(t)dW(t), \quad t\geq0,\\[8pt]
X(0)&=x\in\mathbb{R}^n.\nonumber
\end{align}
In this case, an alternative and simpler proof than the one already presented can be given, where
we show that the mapping $\Gamma_1$ now is a contraction under the norm $||\cdot||_{\lambda}$ for any $\lambda$ satisfying
(\ref{INEQL}). Indeed, for given processes $\bar{X}_i(\cdot)\in\mathbb{M}_{\lambda}^2(0,\infty;\mathbb{R}^n), \, i=1,2,$ we deploy a reasoning similar to the one followed in the proof of Theorem \ref{thm:exuniqgeneral} for the unique solutions $(Y_i,Z_{i})$ and $X_{i}(\cdot), \, i=1,2,$ of the corresponding BSDE (\ref{BSDE}) and the simplified FSDE
\begin{equation*}
\begin{split}
dX(t)&= b\big(t,X(t), Y(t), Z(t)\big) dt-\sigma\big(t,X(t), Y(t)\big)dW(t), \quad t> 0,\\[8pt]
X(0)&=x\in\mathbb{R},
\end{split}
\end{equation*}
according to Lemmata \ref{lem:existBSDE} and \ref{lem:existFSDE}, respectively. Then keeping the notation of (\ref{defabc}), we arrive at the relationships

\begin{align}\label{ineqbarabcsansz}
a\leq \frac{1}{\bar{\lambda}_1}\big[(k_1\epsilon_{1}^{-1}+k_4)b+&k_2\epsilon_2^{-1}c\big],\qquad b\leq\frac{c_1\rho_1^{-1}}{\bar{\lambda}_2}\,\bar{a},
\qquad c\leq\frac{c_1\rho_1^{-1}} {1-c_2\rho_2^{-1}}\,\bar{a},
\end{align}
which are the analogues of those in (\ref{ineqbarabc}) for $k_5=0$. Putting together these inequalities we have that
\begin{align*}
a\leq \frac{c_1\rho_1^{-1}(k_1\epsilon_1^{-1}+k4)}{\bar{\lambda}_{1}}\left[\frac{1}{\bar{\lambda}_{2}}+\frac{k_2\epsilon_2^{-1}}
{(k_1\epsilon_1^{-1}+k_4)(1-c_2\rho_2^{-1})}\right]\bar{a},
\end{align*}
\vskip5pt
\noindent
which shows that the mapping $\Gamma_1$ is a contraction under (\ref{INEQL}) and the same assumptions on the constant parameters. %and the proof is complete.\hfill $\square$
\end{remark}
\begin{remark}\label{rem:lambdarho0}
Due to Proposition \ref{Prop:FBSDEsToREM}, under Hypothesis \ref{hyprem} the stochastic saddlepoint system of (\ref{REM}) is equivalent to
the infinite horizon FBSDEs system of (\ref{FBSDEREM}) and thanks to Theorem  \ref{thm:exuniqgeneral} and Remark \ref{rem:fbsderem} obtains a unique adapted solution, in the sense of Definition \ref{DEF:SOLREM},
when $\lambda$ satisfies (\ref{INEQL}) that simplifies to
\begin{equation*}
2\rho_0+2\mu_0+c_0\rho_1+c_0\rho_1^{-1}<\lambda<-2\mu_1-k_3-k_1\epsilon_1-k_1\epsilon_1^{-1}-k_4.
\end{equation*}
Hence, this relationship formulates explicitly the lower and upper bounds subject to which the stochastic saddlepoint system of (\ref{REM})
is well-posed, according to its solvability as imposed by Definition \ref{DEF:SOLREM}. Furthermore, the first of this double inequality shows exactly how the adjusting rate $\lambda$ should be related with the universal bound $\rho_0$ of the stochastic discount factor
$r(\cdot)$, which is employed to the future expectations model (\ref{REM}) and is needed in the economy for getting to equilibrium via the existence and uniqueness of
a solution for the infinite horizon system of FBSDEs (\ref{FBSDEREM}).
\end{remark}

We conclude this subsection with an example that illustrates the use of Theorem \ref{thm:exuniqgeneral} on the well-posedness of the Dornbusch model presented in Subsection \ref{DORNBUSCH-SECTION}. More precisely, explicit conditions on the coefficients of the Dornbusch saddlepoint system are provided, under which the equilibrium states are eventually achieved. Similar conditions may be given for the rest economic models of Section \ref{sec:rem}, as well.

\begin{example}
As an application of Theorem \ref{thm:exuniqgeneral}, the Dornbusch saddlepoint system that is equivalent to the infinite horizon FBSDEs of (\ref{DORNBUSCH-FBSDE}) achieves a unique adapted solution $(p,s,z)\in\mathbb{M}^2_\lambda(0,\infty;\mathbb{R}^3)$ which drives the system to equilibrium, as long as there exists a constant $\epsilon_{1}>0$ such that $\lambda$ satisfies the condition
%\vskip2pt
\begin{align*}
-\frac{2\nu\eta}{D}+\frac{2|1-\nu\eta-\phi\vartheta|}{|D|}< \lambda < \frac{2\phi(\vartheta+\xi\eta)}{D}-\frac{|\phi\xi|\eta}{|D|}\epsilon_{1}-\frac{|\phi\xi|\eta}{|D|}\epsilon_{1}^{-1}.
\end{align*}
\end{example}

\vskip15pt
\subsection{Comparison of solutions and their continuous dependence on a parameter}
 All the up-to-date comparison results for adapted solutions of FBSDEs in the non Markovian setting (cf. Wu \cite{Wu99},
Wu and Xu \cite{Wu-Xu},
Peng and Shi \cite{Peng-Shi},
and Zhang and Shi \cite{Zhang-Shi})
where obtained via duality techniques and subject to certain traditional monotonicity conditions on
the coefficients of the system. Under our simple and natural monotonicity conditions (H2) of Hypothesis \ref{Hyp},
we provide in the sequel a comparison result for the infinite horizon FBSDEs system of (\ref{GFBSDEs}) with one-dimensional backward process at initial time. This result will
allow us in the next section to derive a comparison property for the introduced notion of a stochastic viscosity solution to the BSPDE that is related
via the four-step-scheme to the previous FBSDEs system%via the four-step-scheme methodology as presented in Subsection \ref{into:subfssrc}
; see Remark \ref{rem:compar-viscos}.

\begin{theorem}\label{thm:comparison}
Let $m=1$ and consider the two infinite horizon fully coupled systems of FBSDEs:
\begin{align}\label{GFBSDEsi}
dX_i(t)&=b_i\big(t,X_i(t),Y_i(t),Z_i(t)\big)dt+\sigma\big(t,X_i(t),Y_i(t),Z_i(t)\big)dW(t), \quad t> 0,\nonumber\\[8pt]
dY_i(t)&=-f_i\big(t,X_i(t),Y_i(t),Z_i(t)\big)dt+Z_i(t)dW(t), \quad t\geq 0,\\[8pt]
X_i(0)&=x\in\mathbb{R}^{n}\nonumber
\end{align}
%corresponding to $X^1(0)=x_1\in\mathbb{R}^{n},$ and  $X^2(0)=x_2\in\mathbb{R}^{n},$ respectively.
for  $i=1,2$. If  the random fields $b_i,$ $f_i,$ for $i=1,2,$ and $\sigma$ satisfy Hypothesis \ref{Hyp}
and
\begin{equation}\label{eqn:bifisigma}
%\begin{split}
b_1(\cdot,x,y,z)\geq b_2(\cdot,x,y,z), \qquad
f_1(\cdot,x,y,z)\geq f_2(\cdot,x,y,z)
%&\sigma^1=\sigma^2=\sigma,
%\end{split}
\end{equation}
hold for every $(x,y,z)\in\mathbb{R}^n\times\mathbb{R}\times \mathbb{R}^{1\times d}$
 then $Y_1(0)\geq Y_2(0)$.
\end{theorem}
\smallskip
\noindent
\textbf{Proof.} According to Theorem \ref{thm:exuniqgeneral} the system of FBSDEs (\ref{GFBSDEsi}) has a unique adapted solution $(X_i, Y_i, Z_i)$ for
$i=1,2,$ respectively, and we define the processes
\begin{align*}
\hat{X}(\cdot)&=X_2(\cdot)-X_1(\cdot), \quad \hat{Y}(\cdot)=Y_2(\cdot)-Y_1(\cdot), \quad \hat{Z}=Z_2(\cdot)-Z_1(\cdot).
%\hat{b}&=b_2-b_1, \qquad\qquad \hat{f}=f_2-f_1.
\end{align*} Applying It\^{o}'s formula on $[0, T]$ to $e^{\lambda\, t}\big[\big(\hat{Y}(t)\big)^{+}\big]^2,$ $t\geq 0,$ taking expectations, letting $T\rightarrow \infty $
and making use of Fatou's lemma, we are led to
\begin{align*}
E\Big[\big(\hat{Y}(0)\big)^{+}\Big]^2\leq & \ (-\lambda+2\mu_2+c_1\rho_1+c_2\rho_2)E\int_{0}^{\infty}e^{\lambda\, s}\ \textbf{1}_{\{\hat{Y}(s)>0\}}\big[\hat{Y}(s)\big]^2ds\\[10pt]
&+(c_2\rho_2^{-1}-1)E\int_{0}^{\infty}e^{\lambda\, s} \ \textbf{1}_{\{\hat{Y}(s)>0\}}||\hat{Z}(s)||^2ds\\[10pt]
&+c_1\rho_1^{-1}E\int_{0}^{\infty}e^{\lambda\, s} \ \textbf{1}_{\{\hat{Y}(s)>0\}}\big|\hat{X}(s)\big|^2ds,
%&-2\int_{t}^{\infty}e^{\lambda\, s}|\big(\hat{Z}^{j}(s)\big)^{+}|^2dW(s)
%&-\int_{t}^{\infty}e^{\lambda\, s}|\hat{Z}^{j}(s)|^2 \mathbb{1}_{\hat{Y}^{j}(s)>0}dW(s)
\end{align*}
\vskip3pt
\noindent
where we have also used an elementary algebraic inequality,
the second relationship of  (\ref{eqn:bifisigma}) and  (H1)-(H3) of Hypothesis \ref{Hyp} for the random fields  $b_i,$ $f_i,$ $i=1,2.$
Following a similar reasoning for $e^{\lambda\, t}\big|\hat{X}(t)\big|^2$, $t\geq 0$, we may also obtain that
%on $[0,T]$ and employing Fatou's lemma  we end up to the relationship
\begin{align*}
E\int_{0}^{\infty}e^{\lambda\, s}|\hat{X}(s)|^2ds\leq \frac{k_1 \epsilon_1^{-1}}{\bar{\lambda}_{1}^{\ast}} E\int_{0}^{\infty}e^{\lambda\, s}\big[\hat{Y}(s)\big]^2ds
+\frac{k_2 \epsilon_2^{-1}}{\bar{\lambda}_{1}^{\ast}} E\int_{0}^{\infty}e^{\lambda\, s}||\hat{Z}(s)||^2ds,
\end{align*}
\vskip3pt
\noindent
where $\bar{\lambda}_{1}^{\ast}\triangleq-\lambda-2\mu_1-k_1\epsilon_1-k_2\epsilon_2>0$ from (\ref{INEQL}).

 Coupling the above inequalities we finally have that
\begin{align*}
E\Big[\big(\hat{Y}(0)\big)^{+}\Big]^2\leq & \left(-\lambda+2\mu_2+c_1\rho_1+c_2\rho_2+c_1\rho_1^{-1}\frac{k_1\epsilon_1^{-1}}{\bar{\lambda}_{1}^{\ast}}\right)||\hat{Y}||_{\lambda}^2\\[13pt]
& +\left(c_2\rho_2^{-1}+c_1\rho_1^{-1}\frac{k_2\epsilon_2^{-1}}{\bar{\lambda}_{1}^{\ast}}-1\right)||\hat{Z}||_{\lambda}^2,
\end{align*}
\vskip3pt
\noindent
where from (H5) of Hypothesis \ref{Hyp} we can easily see that
\begin{align*}
-\lambda+2\mu_2+c_1\rho_1+c_2\rho_2+c_1\rho_1^{-1}\frac{k_1\epsilon_1^{-1}}{\bar{\lambda}_{1}^{\ast}}<0,\qquad
 c_2\rho_2^{-1}+c_1\rho_1^{-1}\frac{k_2\epsilon_2^{-1}}{\bar{\lambda}_{1}^{\ast}}-1<0 .
\end{align*}
%\vskip5pt
%\noindent
Therefore we get immediately that $\big(\hat{Y}(0)\big)^{+}=0,$ which concludes the proof. \hfill $\square$
\smallskip

We shall investigate next the continuous dependence of solutions of such infinite horizon FBSDEs on a parameter. If this is the case, it follows immediately that these FBSDEs are well-posed in the sense of J. Hadamard,
since they not only admit a unique solution, thanks to Theorem \ref{thm:exuniqgeneral}, but its behavior also changes continuously with the initial conditions. This results in a better understanding of how
sensitive are the assets $Y(\cdot)$ to the fundamentals $X(\cdot)$ in our economic model; in other words, the different trajectories of the
fundamentals designate in a smooth manner the corresponding trajectories of the assets. Furthermore, in the next section we shall make use of this continuity property to establish the appropriate regularity for a suggested stochastic viscosity solution of the associated BSPDE; cf. Proposition \ref{thm:svs}.

 Assume that Hypothesis \ref{Hyp} holds for a family $\big\{x(\alpha), b(\alpha\,;\cdot\,), \sigma(\alpha\,;\cdot\,), f(\alpha\,; \cdot\,)$, $\alpha\in\mathbb{R}\big\}$ of initial conditions and random coefficients of the infinite horizon FBSDEs system (\ref{GFBSDEs}), uniformly in the parameter $\alpha.$ Furthermore, Theorem \ref{thm:exuniqgeneral} establishes a unique adapted solution triplet for each one of these systems, which we denote accordingly by $(X_\alpha, Y_\alpha, Z_\alpha).$

\begin{theorem} \label{thm:cdsp}
If the family of the initial conditions and the random coefficients of the infinite horizon FBSDES (\ref{GFBSDEs}) is Lipschitz continuous with respect to the parameter $\alpha$, i.e.,
\begin{align*}
|x(\alpha_2)-x(\alpha_1)|&\leq C_1|\alpha_2-\alpha_1|,\nonumber\\[15pt]
\big|\big|b\big(\alpha_2\,;\cdot\,,X_{\alpha_1}(\cdot),Y_{\alpha_1}(\cdot),Z_{\alpha_1}(\cdot)\big)-b\big(\alpha_1\,;\cdot\,,X_{\alpha_1}(\cdot),Y_{\alpha_1}(\cdot),Z_{\alpha_1}(\cdot)\big)\big|\big|_\lambda
&\leq C_1|\alpha_2-\alpha_1|,\\[15pt]
\big|\big|\sigma\big(\alpha_2\,;\cdot\,,X_{\alpha_1}(\cdot),Y_{\alpha_1}(\cdot),Z_{\alpha_1}(\cdot)
\big)-\sigma\big(\alpha_1\,;\cdot\,,X_{\alpha_1}(\cdot),Y_{\alpha_1}(\cdot),Z_{\alpha_1}(\cdot)
\big)\big|\big|_\lambda&\leq C_1|\alpha_2-\alpha_1|,\\[15pt]
\big|\big|f\big(\alpha_2\,;\cdot\,,X_{\alpha_1}(\cdot),Y_{\alpha_1}(\cdot),Z_{\alpha_1}(\cdot)\big)-f\big(\alpha_1\,;\cdot\,,X_{\alpha_1}(\cdot),Y_{\alpha_1}(\cdot),Z_{\alpha_1}(\cdot)\big)\big|\big|_\lambda&\leq C_1|\alpha_2-\alpha_1|\nonumber
\end{align*}
\vskip5pt
\noindent
for every $\alpha_1,$ $\alpha_2\in\mathbb{R}$ and some constant $C_1>0$, then we have that
\begin{align*}
 E \left(\sup_{t\geq 0}|X_{a_2}(t)-X_{a_1}(t)|^2\right)&\leq C_2|a_2-a_1|^2, \\[10pt]
 E \left(\sup_{t\geq 0}|Y_{a_2}(t)-Y_{a_1}(t)|^2\right)&\leq C_2|a_2-a_1|^2,\\[10pt]
 ||Z_{a_2}-Z_{a_1}||_\lambda &\leq C_2|a_2-a_1|
\end{align*}
for some other constant $C_2>0$.
\end{theorem}
\vskip5pt
\noindent
\textbf{Proof.}
Given any $t\geq 0$, $(x,y,z)\in\mathbb{R}^n\times\mathbb{R}^m\times\mathbb{R}^{m\times d}$ and $\alpha_1,$ $\alpha_2\in\mathbb{R}$, we define the processes
 \[ \delta X_{\alpha}(\cdot)\triangleq X_{\alpha_2}(\cdot)-X_{\alpha_1}(\cdot), \quad \delta Y_{\alpha}(\cdot)\triangleq  Y_{\alpha_2}(\cdot)-Y_{\alpha_1}(\cdot),  \quad \delta Z_{ \alpha}(\cdot)\triangleq Z_{\alpha_2}(\cdot)-Z_{\alpha_1}(\cdot),\]
the random fields

\begin{align*}
\delta B_{\alpha}(t,x,y,z)&\triangleq b\big(\alpha_{2};t,x+X_{\alpha_1}(t),y+Y_{\alpha_1}(t),z+Z_{\alpha_1}(t)\big)\\
&\qquad\qquad \qquad\qquad \qquad\, -b\big(\alpha_{1};t,X_{\alpha_1}(t),Y_{\alpha_1}(t),Z_{\alpha_1}(t)\big),\\[10pt]
\delta\Sigma_{\alpha}(t,x,y,z)&\triangleq \sigma\big(\alpha_{2};t,x+X_{\alpha_1}(t),y+Y_{\alpha_1}(t),z+Z_{\alpha_1}(t)\big)\\
&\qquad\qquad \qquad\qquad \qquad\, -\sigma\big(\alpha_{1};t,X_{\alpha_1}(t),Y_{\alpha_1}(t),Z_{\alpha_1}(t)\big),\\[10pt]
\delta F_{\alpha}(t,x,y,z)&\triangleq f\big(\alpha_{2};t,x+X_{\alpha_1}(t),y+Y_{\alpha_1}(t),z+Z_{\alpha_1}(t)\big)\\
&\qquad\qquad \qquad\qquad \qquad\, -f\big(\alpha_{1};t,X_{\alpha_1}(t),Y_{\alpha_1}(t),Z_{\alpha_1}(t)\big),
\end{align*}
and the function
\[
\delta x(\alpha)\triangleq x(\alpha_2)-x(\alpha_1).\]
Then the process triplet  $(\delta X_\alpha,\delta Y_\alpha, \delta Z_\alpha)$ solves the system of infinite horizon random FBSDEs
\begin{align}\label{FBSDEdiff}
dX(t)&= \delta B_\alpha\big(t,X(t),Y(t),Z(t)\big)dt+ \delta\Sigma_\alpha\big(t,X(t),Y(t),Z(t)\big)dW(t),\quad t>0,\nonumber\\[5pt]
dY(t)&=- \delta F_\alpha\big(t,X(t),Y(t),Z(t)\big)dt+ Z(t) dW(t),\quad t\geq 0,\\[5pt]
X(0)&=\delta x(\alpha).\nonumber
\end{align}
We devote the rest of the proof to the derivation of the estimate
\begin{align}\label{metrxyz}
||\delta X_{\alpha}||_{\lambda}^2+||\delta Y_{\alpha}||_{\lambda}^2+||\delta Z_{ \alpha}||_{\lambda}^2\leq C_3\bigg[&|\delta x(\alpha)|^2+||\delta B_{\alpha}(\cdot,0,0,0)||_{\lambda}^2\\[5pt]
&+||\delta F_{\alpha}(\cdot,0,0,0)||_{\lambda}^2+||\delta\Sigma_{\alpha}(\cdot,0,0,0)||_{\lambda}^2\bigg],\nonumber
\end{align}
where $\lambda$ satisfies (\ref{INEQL}) and $C_3>0$ is a constant which depends on $k_{i}, \ i=1,...,4,$ and $c_i,\ \mu_i,$ $i=1,2;$
the assertion of the theorem follows then immediately from its hypotheses, Burkholder's inequality and Fatou's lemma.

Since the process triplet $(\delta X_\alpha,\delta Y_\alpha, \delta Z_\alpha)$ satisfies the FSDE of (\ref{FBSDEdiff}), a combination of It\^{o}'s lemma with an elementary algebraic inequality,
 the inequality (\ref{INEQL}) and Hypothesis \ref{Hyp} for the random fields $b$, $\sigma$ leads for all $\epsilon, T>0$ to
\begin{align*}%\label{ineqx1}
\bar{\lambda}_1 E\int_{0}^{T} e^{\lambda t}|\delta X_{\alpha}(t)|^2dt\leq \ & |\delta x(\alpha)|^2+\big[k_1\epsilon_{1}^{-1}+k_4(1+\epsilon)\big]E\int_{0}^{T} e^{\lambda t}|\delta Y_{\alpha}(t)|^2dt\nonumber\\[10pt]
&+\big[k_2\epsilon_2^{-1}+k_5(1+\epsilon)\big]E\int_{0}^{T} e^{\lambda t}||\delta Z_{ \alpha}(t)||^2dt\nonumber\\[10pt]
&+\epsilon^{-1}E\int_{0}^{T} e^{\lambda t}|\delta B_{\alpha}(t,0,0,0)|^2dt\\[10pt]
&+(1+\epsilon^{-1})E\int_{0}^{T} e^{\lambda t}||\delta\Sigma_{\alpha}(t,0,0,0)||^2dt\nonumber
,\nonumber
\end{align*}
 where we have also assumed that
 \begin{equation}\label{barl11}
 \bar{\lambda}_1\triangleq -\lambda-2\mu_1-k_1\epsilon_1-k_2\epsilon_2-k_3(1+\epsilon)-\epsilon>0.
\end{equation}
Then Fatou's lemma implies the inequality
\begin{align}\label{metrx1}
||\delta X_{\alpha}||_{\lambda}^2\leq \frac{1}{\bar{\lambda}_1}\bigg\{&\big[k_1\epsilon_{1}^{-1}+k_4(1+\epsilon)\big]||\delta Y_{\alpha}||_{\lambda}^2
+\big[k_2\epsilon_2^{-1}+k_5(1+\epsilon)\big]
||\delta Z_{\alpha}||_{\lambda}^2\\[10pt]
&+|\delta x(\alpha)|^2+\epsilon^{-1}||\delta B_{\alpha}(\cdot,0,0,0)||_{\lambda}^2
+(1+\epsilon^{-1})||\delta \Sigma_{\alpha}(\cdot,0,0,0)||_{\lambda}^2\bigg\}.\nonumber
\end{align}

On the other hand, $(\delta X_{\alpha}, \delta Y_{\alpha},\delta Z_{\alpha})$ satisfies the BSDE of (\ref{FBSDEdiff}) as well. Henceforth, following the same reasoning as above
and assuming that
\begin{equation}\label{barl21}
 \bar{\lambda}_2\triangleq \lambda-2\mu_2-c_1\rho_1-c_2\rho_2-\epsilon>0,
\end{equation}
we end for any $t\geq 0$ in

\begin{align*}
&\quad \ \ e^{\lambda t}E|\delta Y_{\alpha}(t)|^2+\bar{\lambda}_2 E\int_{t}^{\infty} e^{\lambda s}|\delta Y_{\alpha}(s)|^2ds+(1-c_2\rho_{2}^{-1})E\int_{t}^{\infty} e^{\lambda s}||\delta Z_{\alpha}(s)||^2ds\\[15pt]
&\leq c_1\rho_{1}^{-1}E\int_{t}^{\infty} e^{\lambda s}|\delta X_{\alpha}(s)|^2ds +\epsilon^{-1}E\int_{t}^{\infty} e^{\lambda s}|\delta F_\alpha(s,0,0,0)|^2ds.
\end{align*}
\vskip3pt
\noindent
From the above inequality we immediately obtain that
\begin{align}\label{metryz1}
\begin{split}
||\delta Y_{\alpha}||_{\lambda}^2&\leq \frac{1}{\bar{\lambda}_2}\bigg[c_1\rho_{1}^{-1}||\delta X_{\alpha}||_{\lambda}^2+\epsilon^{-1}||\delta F_\alpha(\cdot,0,0,0)||_{\lambda}^2\bigg],\\[15pt]
||\delta Z_{\alpha}||_{\lambda}^2&\leq \frac{1}{{1-c_2\rho_2^{-1}}}\bigg[c_1\rho_{1}^{-1}||\delta X_{\alpha}||_{\lambda}^2+\epsilon^{-1}||\delta F_\alpha(\cdot,0,0,0)||_{\lambda}^2\bigg].
\end{split}
\end{align}
\smallskip

We ease notation by setting
\begin{align*}
a&\triangleq ||\delta X_{\alpha}||_{\lambda}^2,\qquad \quad b\triangleq ||\delta Y_{\alpha}||_{\lambda}^2,\qquad \quad c\triangleq ||\delta Z_{\alpha}||_{\lambda}^2, \nonumber\\[10pt]
a_0&\triangleq |\delta x(\alpha)|^2+\epsilon^{-1}||\delta B_{\alpha}(\cdot,0,0,0)||_{\lambda}^2
+(1+\epsilon^{-1})||\delta \Sigma_{\alpha}(\cdot,0,0,0)||_{\lambda}^2,\\[15pt]
b_0&\triangleq \epsilon^{-1}||\delta F_\alpha(\cdot,0,0,0)||_{\lambda}^2,\nonumber
\end{align*}
and reformulate the inequalities (\ref{metrx1}) and (\ref{metryz1}), respectively, as
\begin{align} \label{ineqabc1}
\begin{split}
a&\leq \frac{1}{\bar{\lambda}_1}\bigg\{\big[k_1\epsilon_{1}^{-1}+k_4(1+\epsilon)\big]b
+\big[k_2\epsilon_2^{-1}+k_5(1+\epsilon)\big]
c +a_0\bigg\}\\[10pt]
b&\leq \frac{1}{\bar{\lambda}_2}\big[c_1\rho_{1}^{-1}a+b_0\big],\qquad
c\leq \frac{1}{{1-c_2\rho_2^{-1}}}\big[c_1\rho_{1}^{-1}a+b_0\big].
\end{split}
\end{align}
\vskip5pt
\noindent
Making use of these inequalities and denoting $\gamma^{\ast}\triangleq\frac{k_2\epsilon_2^{-1}+k_5(1+\epsilon)}{k_1\epsilon_1^{-1}+k_4(1+\epsilon)}$ we get that

\begin{align}\label{INEQbcgammast}
b+\gamma^{\ast} c\leq & \ \frac{c_1\rho_1^{-1}}{\bar{\lambda}_1}\left[\frac{1}{\bar{\lambda}_2}+\frac{\gamma^{\ast}}{1-c_2\rho_2^{-1}}\right] \big[k_1\epsilon^{-1}+k_4(1+\epsilon)\big](b+\gamma^{\ast} c) \nonumber\\\\
&+\left[\frac{1}{\bar{\lambda}_2}+\frac{\gamma^{\ast}}{1-c_2\rho_2^{-1}}\right]
\left[\frac{c_1\rho_1^{-1}a_0}{\bar{\lambda}_1}+b_0\right]\nonumber.
\end{align}
\vskip5pt
\noindent
 In view of (\ref{INEQcrho}) and (\ref{INEQL}),  we may choose a sufficiently small $\epsilon>0$ such that the inequalities in (\ref{barl11}), (\ref{barl21}) and
 \[\mu\triangleq \frac{c_1\rho_1^{-1}}{\bar{\lambda}_1}\left[\frac{1}{\bar{\lambda}_2}+\frac{\gamma^{\ast}}{1-c_2\rho_2^{-1}}\right] \big[k_1\epsilon^{-1}+k_4(1+\epsilon)\big]<1\]
\vskip5pt
\noindent
hold. Therefore, the inequality (\ref{INEQbcgammast}) transforms to
\[b+\gamma^{\ast} c\leq \frac{1}{1-\mu}\left[\frac{1}{\bar{\lambda}_2}+\frac{\gamma^{\ast}}{1-c_2\rho_2^{-1}}\right]  \left[\frac{c_1\rho_1^{-1}a_0}{\bar{\lambda}_1}+b_0\right],\]
\vskip5pt
\noindent
which together with the first inequality of (\ref{ineqabc1}) justify (\ref{metrxyz}).
\hfill $\square$
\smallskip
%\eject

This result leads immediately to the following corollary, according to which the adapted solution triplet $(X, Y, Z)$ of the infinite horizon random FBSDEs (\ref{GFBSDEs}) depends continuously on its initial conditions $(t,x).$
\begin{corollary}\label{col:cdsp}
Suppose that Hypothesis \ref{Hyp} holds. If $(X_{t_i,x_i}, Y_{t_i,x_i}, Z_{t_i,x_i})$ is the adapted solution of (\ref{GFBSDEs}) associated to the initial
point $(t_i,x_i), i=1,2,$ (cf. (\ref{eqn:FBSDEsanss}) below) then we get

\begin{equation}
\begin{aligned}\label{COLcpd}
&\quad E\left( \sup_{s\geq t_1\vee t_2 }e^{\lambda s} |X_{t_1,x_1}(s)-X_{t_2,x_2}(s)|^2\right)+E \left(\sup_{s \geq t_1\vee t_2}e^{\lambda s} |Y_{t_1,x_1}(s)-Y_{t_2,x_2}(s)|^2\right) \\[20pt]
&\quad+E\left(\int_{t_1\vee t_2}^{\infty}  e^{\lambda s} |Z_{t_1,x_1}(s)-Z_{t_2,x_2}(s)|^2ds\right) \leq C \bigg(|x_1-x_2|^2+(1+|x_1|^2\vee|x_2|^2)|t_1-t_2|\bigg).
\end{aligned}
\end{equation}
\end{corollary}

\vskip 20pt
\section{(In-)Finite Horizon Random FBSDEs And (Stationary)\\ Stochastic Viscosity Solutions}\label{sub:svs}
It is well known that in a finite time-horizon Markovian framework one may employ the theory of BSDEs to obtain a probabilistic representation for the solution of a linear or even semi-linear second order parabolic PDE
by means of the Feynman-Kac formula; cf. Peng \cite{Peng-91}, \cite{Peng-92}
and Pardoux et al. \cite{Pardoux-97}.
By the same token, a fully coupled system of FBSDEs is related to a quasilinear parabolic backward PDE, as argued e.g. in Ma et al. \cite{Ma-Protter-Yong}
and Pardoux and Tang \cite{Pardoux-Tang},
while a second-order BSDE (2BSDE) is associated to a fully nonlinear PDE, according to
Cheridito et al. \cite{Cheriditoetal}
and Soner et al. \cite{Soneratal}.
In this vein, Pardoux and Tang \cite{Pardoux-Tang}
took advantage of the theory of FBSDEs to provide a probabilistic characterization for the solution of a quasilinear PDE of parabolic type, in the \emph{viscosity} sense.

In this section, we consider the \emph{random} FBSDEs system of (\ref{FBSDEssigma}) for an arbitrary finite time-horizon
and by making use of the semimartingale decomposition of a random field we
introduce a new notion, to the best of our knowledge, of a \emph{stochastic viscosity solution} for the corresponding quasilinear BSPDE of parabolic type.
Motivated further by Zhang and Zhao \cite{Zhang-Zhao},
%\cite{Zhang-Zhao}
who studied the connection of a stationary solution for an SPDE with the solution of a backward doubly SDEs system,
we provide additional assumptions on the coefficients of the random FBSDEs system (\ref{FBSDEssigma}) under which its solution becomes
\emph{stationary} and leads in turn to the existence of a stochastic viscosity solution for
the associated BSPDE with an indistinguishable version that is \emph{stationary} as well.

We shall consider first the simpler case of $m=1$. Let Hypothesis \ref{Hyp} hold with $k_5=0,$ that is the random field $\sigma$ is independent of the parameter $z$, and for every initial condition $(t,x)\in \mathbb{R}_{+}\times \mathbb{R}^{n}$,
according to Remark \ref{rem:exisuniqsigma}, denote by
$(X_{t,x}, Y_{t,x}, Z_{t,x})$ the unique adapted solution of class $\mathbb{M}_{\lambda}^2(t,\infty;\mathbb{R}^n\times\mathbb{R}\times\mathbb{R}^{1\times d})$ to the infinite horizon system of random FBSDEs
\begin{align}\label{eqn:FBSDEsanss}
dX_{t,x}(s)&=b\big(s,X_{t,x}(s),Y_{t,x}(s), Z_{t,x}(s)\big)ds+\sigma\big(s,X_{t,x}(s), Y_{t,x}(s)\big)dW(s),\quad s>t,\nonumber\\[10pt]
dY_{t,x}(s)&=-f\big(s,X_{t,x}(s), Y_{t,x}(s), Z_{t,x}(s)\big)ds+ Z_{t,x}(s)dW(s), \quad s\geq t,  \\[10pt]
X_{t,x}(t)&=x,\nonumber
\end{align}
and its finite horizon counterpart

\begin{align}\label{eqn:FBSDEsanss-FINITE}
dX_{t,x}(s)&=b\big(s,X_{t,x}(s),Y_{t,x}(s), Z_{t,x}(s)\big)ds+\sigma\big(s,X_{t,x}(s), Y_{t,x}(s)\big)dW(s),\quad t<s\leq T,\nonumber\\[10pt]
dY_{t,x}(s)&=-f\big(s,X_{t,x}(s), Y_{t,x}(s), Z_{t,x}(s)\big)ds+ Z_{t,x}(s)dW(s), \quad t\leq s<T,  \\[10pt]
X_{t,x}(t)&=x, \,\,\, Y_{t,x}(T)=\Phi(T,X_{t,x}(T)),\nonumber
\end{align}
where $\Phi$ is a given random field that connects the forward and the backward variables at time $T$ and
the solution   $(X_{t,x},Y_{t,x},Z_{t,x})$ belongs now to $ \mathbb{M}_{\lambda}^2(t,T;\mathbb{R}^n\times\mathbb{R}\times\mathbb{R}^{1\times d})$.

\subsection{The finite horizon case}

We first consider the finite horizon system  \eqref{eqn:FBSDEsanss-FINITE} for a given (arbitrary) random field $\Phi$.
Following the four-step-scheme methodology %first developed by Ma et al. \cite{Ma-Protter-Yong} and later generalized
 for random fields by Ma and Yong \cite{Ma-Yong97},
  it would be interesting to have a random field $v$ that connects the instantaneous state of the variable $X(\cdot)$ with the instantaneous state of the variable $Y(\cdot)$, and accordingly with the variable $Z(\cdot)$, i.e.
such that $Y(s)=v\big(s,X(s)\big)$ for every $t\leq s\leq T$. This \emph{representing random field} should be very useful for the solution of the FBSDEs \eqref{eqn:FBSDEsanss-FINITE}, as it would allow the decoupling of the system and eventually lead to the treatment of the corresponding BSPDE.

If the above case is not possible, the next to do would be to look instead for a random field $v'$ that connects the initial value $x$ of the forward variable $X(\cdot)$, which is known, with the initial value $y$
of the backward variable $Y(\cdot)$, which is unknown and can only be obtained after we solve the FBSDEs problem. This random field would be such that $y=v'(t,x)$ and the knowledge of $y$ would allow us to solve the system using a shooting algorithm, in which we would treat the process pair $(X,Y)$ as the solution of a forward system with initial condition $(x,y)$. This random field $v'$ is not a true representing random field in the sense that it holds only for the initial values and not for the instantaneous values as in the previous case.

It would be shown shortly (cf. Theorem \ref{thm:class-visc}) that the random field we seek in each case is related to the solution of the following quasilinear second-order parabolic BSPDE:
\begin{align}\label{equ:BSPDEu}
v(t,x)&=\Phi(T,x)+\int_{t}^{T}\bigg[(Lv)\Big(s,x,v(s,x),\left(\nabla_x v\right)^{\top}(s,x)\,\sigma\big(s,x,v(s,x)\big)+\Psi^v(s,x)\Big)\nonumber\\[15pt]
& \qquad\qquad\qquad\quad +f\Big(s,x,v(s,x),\left(\nabla_x v\right)^{\top}(s,x)\,\sigma\big(s,x,v(s,x)\big)+\Psi^{v}(s,x)\Big)\\[15pt]
& \qquad\qquad\qquad\quad +tr\left\{\left(\nabla_x \Psi^{v}\right)^{\top}(s,x)\,\sigma\big(s,x,v(s,x)\big)\right\}\bigg]ds -\int_{t}^{T}\Psi^{v}(s,x)dW(s),\nonumber
\end{align}
which holds for every $(t,x)\in [0,T]\times\mathbb{R}^{n}$.  In the above, we have used the random elliptic operator $L$ defined for every  $(t,\bar{x},y,z)\in [0,T)\times\mathbb{R}^{n}\times\mathbb{R}\times\mathbb{R}^{1\times d}$  by
\begin{equation*}
(Lv)(t,\bar{x},y,z)\triangleq \frac{1}{2}\sum_{i,j=1}^{n} \big((\sigma\sigma^{\top})(t,\bar{x},y)\big)_{i,j}\, \frac{\partial^2v}{\partial x_i \partial x_j} (t,\bar{x})+\big\langle b(t,\bar{x},y,z),\nabla_x v(t,\bar{x})\big\rangle.
\end{equation*}
We are also going to consider the following two types of solutions.

\begin{definition}\label{def:svs}
(i) The random field pair $(v, \Psi^v)$ is a \emph{classical solution} of the BSPDE (\ref{equ:BSPDEu}) if it is of class $C_{\mathbb{F}}\big([0, T];\mathbb{L}^2(\Omega;C^3(\mathbb{R}^{n};\mathbb{R}))\big) \times \mathbb{M}_{\lambda}^2\big([0, T];C^2(\mathbb{R}^{n};\mathbb{R}^{1\times d})\big)$ and satisfies the relation of (\ref{equ:BSPDEu}).

(ii) Let $v$ be a random field of class
 $C_{\mathbb{F}}\big([0,T];\mathbb{L}^2(\Omega;C(\mathbb{R}^{n};\mathbb{R}))\big)$ and admit a semimartingale decomposition. The random field pair
  $(v,\Psi^{v})$ is called a \emph{stochastic viscosity subsolution} (resp. \emph{supersolution})
  of the BSPDE (\ref{equ:BSPDEu}) if whenever the random field triplet $(\phi, A^{\phi}, \Psi^{\phi})$ %with semimartingale representation
 %$$d\phi(t,x)=A^{\phi}(t,x)dt+\Psi^{\phi}(t,x)dW(t),\qquad \forall\quad (t,x)\in \mathbb{R}_{+}\times\mathbb{R}^{n} $$
    belongs to the class of $C_{\mathbb{F}}\big([0,T];\mathbb{L}^2(\Omega;C^2(\mathbb{R}^{n};\mathbb{R}))\big) \times C_{\mathbb{F}}\big([0,T];\mathbb{L}^2(\Omega;C^1(\mathbb{R}^{n};\mathbb{R}))\big) \times\mathbb{M}_{\lambda}^2\big([0,T];C^2(\mathbb{R}^{n};\mathbb{R}^{1\times d})\big)$ and $\,(t,x)\in[0,T]\times\mathbb{R}^{n}\,$ is a local minimum (resp. maximum) of $\phi-v$, then we have that
 \begin{align*}
 A^{\phi}(t,x)&+(L\phi)\Big(t,x,v(t,x),\left(\nabla_x\phi\right)^{\top}(t,x)\,\sigma\big(t,x,v(t,x)\big)+\Psi^\phi(t,x)\Big)\\[13pt]
 &+f\Big(t,x,v(t,x),\left(\nabla_x\phi\right)^{\top}(t,x)\,\sigma\big(t,x,v(t,x)\big)+\Psi^\phi(t,x)\Big)\\[13pt]
 &+tr\left\{\left(\nabla_x\Psi^{\phi}\right)^{\top}(t,x)\,\sigma\big(t,x,v(t,x)\big)\right\}\geq 0
 \end{align*}
 \vskip5pt
 \noindent
\big[ resp.  \begin{align*}
 A^{\phi}(t,x)&+(L\phi)\Big(t,x,v(t,x),\left(\nabla_x\phi\right)^{\top}(t,x)\,\sigma\big(t,x,v(t,x)\big)+\Psi^\phi(t,x)\Big)\\[13pt]
 &+f\Big(t,x,v(t,x),\left(\nabla_x\phi\right)^{\top}(t,x)\,\sigma\big(t,x,v(t,x)\big)+\Psi^\phi(t,x)\Big)\\[13pt]
 &+tr\left\{\left(\nabla_x\Psi^{\phi}\right)^{\top}(t,x)\,\sigma\big(t,x,v(t,x)\big)\right\}\leq 0\, \ \big].
 \end{align*}
\vskip5pt
 \noindent
 We shall call $(v,\Psi^{v})$ a \emph{stochastic viscosity solution} of the BSPDE (\ref{equ:BSPDEu}) if it is both a stochastic viscosity sub- and super- solution.
\end{definition}

According to Proposition \ref{prop:IKW}, the additional order of smoothness considered in the above definition for the classical solution of BSPDE (\ref{equ:BSPDEu}) will allow us below to employ the generalized
 IKW formula. We should also note here that this definition of a stochastic viscosity solution extends the deterministic one proposed by Pardoux and Tang \cite{Pardoux-Tang}
in the Markovian case.
%We  are ready now to prove the next Theorem which deals with the solvability of the  SPDE (\ref{equ:BSPDEu}) thus of the solvability of the FBSDEs (\ref{FBSDEssigma}).

\begin{theorem}\label{thm:class-visc}
Let the random fields $b,$ $\sigma$ and $f$ be globally continuous and satisfy  Hypothesis \ref{Hyp} with $k_5=0.$

(i) Suppose BSPDE (\ref{equ:BSPDEu}) admits a unique classical solution $(v,\Psi^{v})$. Then the forward $X(\cdot)$ and backward $\big(Y(\cdot),Z(\cdot)\big)$ variables of the FBSDEs system (\ref{eqn:FBSDEsanss-FINITE}) are connected via the relationships:

\begin{equation} \label{FSSC}
\begin{split}
Y_{t,x}(s)&=v\big(s,X_{t,x}(s)\big),\ \quad t\leq s\leq T, \\[10pt]
Z_{t,x}(s)&=\Psi^{v}\big(s,X_{t,x}(s)\big)+\left(\nabla_x v\right)^{\top}\big(s,X_{t,x}(s)\big)\,\sigma\Big(s, X_{t,x}(s), v\big(s,X_{t,x}(s)\big)\Big), \  \ t \leq s \leq T;
\end{split}
\end{equation}
in other words, $v$ is the representing random field.

(ii) Suppose BSPDE (\ref{equ:BSPDEu}) admits a unique viscosity solution $(v,\Psi^{v})$. Then the initial value $x$ of the forward variable $X(\cdot)$ of the FBSDEs system (\ref{eqn:FBSDEsanss-FINITE}) is connected with the (unknown) initial value $y$ of its backward variable $Y(\cdot)$ by the relation:
\[ y=v(t,x). \]
\end{theorem}
\noindent
\textbf{Proof.}
For the first assertion of the theorem, it suffices to apply the generalized IKW formula to the process  $v\big(\cdot,X_{t,x}(\cdot)\big)$ on $[t,T]$ and compare the resulting equation with the backward one of the FBSDEs system (\ref{equ:BSPDEu}). The second assertion of the theorem is a direct consequence of the following proposition, which we believe is of interest on its own right. \hfill$\square$

\begin{proposition}\label{thm:svs}
 Suppose that the random fields $b,$ $\sigma$ and $f$ are globally continuous and satisfy  Hypothesis \ref{Hyp} with $k_5=0.$ Then the random field $v$, defined in terms of the unique adapted solution $(X_{t,x},Y_{t,x},Z_{t,x})$ of  the FBSDEs (\ref{eqn:FBSDEsanss-FINITE}) by
 \begin{equation}\label{eqn:vy}
 v(t,x)\triangleq Y_{t,x}(t),\qquad \forall\ (t,x)\in[0,T]\times \mathbb{R}^n,
 \end{equation}
 is of class  $C_{\mathbb{F}}\big([0,T];\mathbb{L}^2(\Omega;C(\mathbb{R}^{n};\mathbb{R}))\big)$.
 If additionally $v$ admits a semimartingale decomposition then the random field pair $(v,\Psi^{v})$ is a stochastic viscosity solution of the BSPDE (\ref{equ:BSPDEu}).
\end{proposition}
%\begin{remark}
%Even though in Proposition \ref{thm:svs} there is no explicit definition for the random field $\Psi^{v}$,
%we shall argue that it is determined uniquely through the definition of the random field $v$ in (\ref{eqn:vy}). Firstly, the required regularity of $v$ follows directly from Corollary \ref{col:cdsp} and the fact that
% $$E\left(\sup_{s\geq t}e^{\lambda\,s}Y_{t,x}^2(s)\right)<\infty$$
%for any $(t,x)\in \mathbb{R}_{+}\times\mathbb{R}^n,$ in accordance with the relationship (\ref{eqn:sup}).
%Then, since $$E\big(v^2(t,x)\big)<\infty$$
%\end{remark}
\noindent
\textbf{Proof.}
The asserted regularity of $v$ given by (\ref{eqn:vy}) follows directly from Corollary \ref{col:cdsp}.
It remains to prove that the pair  $(v, \Psi^{v})$ is a stochastic viscosity solution of the
BSPDE (\ref{equ:BSPDEu}). In the rest of the proof we shall establish that $(v, \Psi^{v})$ is a stochastic viscosity subsolution. One can use analogous arguments to show that $(v, \Psi^{v})$ is a stochastic viscosity supersolution as well.

Let the random field triplet  $(\phi, A^{\phi}, \Psi^{\phi})$ belongs  to the class $C_{\mathbb{F}}\big([0,T];\mathbb{L}^2(\Omega;C^2(\mathbb{R}^{n};\mathbb{R}))\big) \times C_{\mathbb{F}}\big([0,T];\mathbb{L}^2(\Omega;C^1(\mathbb{R}^{n};\mathbb{R}))\big) \times\mathbb{M}_{\lambda}^2\big([0,T];C^2(\mathbb{R}^{n};\mathbb{R}^{1\times d})\big)$ and $\,(t,x)\in[0,T]\times\mathbb{R}^{n}\,$ be a local minimum of $\phi-v$ , and without loss of generality assume  that  $\phi(t,x)=v(t,x).$
If we assume that
\begin{align}\label{eqn:localmin}
 A^{\phi}(t,x)&+(L\phi)\Big(t,x,v(t,x),\left(\nabla_x\phi\right)^{\top}(t,x)\,\sigma\big(t,x,v(t,x)\big)+\Psi^\phi(t,x)\Big)\nonumber\\[10pt]
 &+f\Big(t,x,v(t,x),\left(\nabla_x\phi\right)^{\top}(t,x)\,\sigma\big(t,x,v(t,x)\big)+\Psi^\phi(t,x)\Big)\\[10pt]
 &+tr\left\{\left(\nabla_x\Psi^{\phi}\right)^{\top}(t,x)\,\sigma\big(s,x,v(t,x)\big)\right\}< 0,\nonumber
 \end{align}
it suffices to reach a contradiction to finish the proof.

Since $\phi-v$ has a local minimum at $(t,x)\in [0,T]\times\mathbb{R}^n$ and $v$ is continuous, there exists $0<\alpha<T-t$ such that for all $(s,y)\in[t,T]\times \mathbb{R}^n$ with
$t\leq s\leq t+\alpha$ and $|x-y|\leq \alpha$ we have  that $ v(s,y)\leq \phi(s,y)$ and
\begin{align}\label{eqn:contrad}
A^{\phi}(s,y)&+(L\phi)\Big(s,y,v(s,y),\left(\nabla_x\phi\right)^{\top}(s,y)\,\sigma\big(s,y,v(s,y)\big)+\Psi^\phi(s,y)\Big)\nonumber\\[10pt]
 &+f\Big(s,y,v(s,y),\left(\nabla_x\phi\right)^{\top}(s,y)\,\sigma\big(s,y,v(s,y)\big)+\Psi^\phi(s,y)\Big)\\[10pt]
 &+tr\left\{\left(\nabla_x\Psi^{\phi}\right)^{\top}(s,y)\,\sigma\big(s,y,v(s,y)\big)\right\}< 0.\nonumber
\end{align}
Furthermore, due to the uniqueness of the solution of the  FBSDEs  (\ref{eqn:FBSDEsanss-FINITE}) %and definition (\ref{}
we derive that
\begin{equation*}
Y_{t,x}(s)=Y_{s,X_{t,x}(s)}(s)=v\big(s, X_{t,x}(s)\big), \qquad t\leq s\leq T.
\end{equation*}
\vskip5pt
\noindent
We  define the stopping time
\begin{equation*}
\tau\triangleq \inf\{s>t: |X_{t,x}(s)-x|\geq \alpha\}\wedge (t+\alpha)
\end{equation*}
\vskip5pt
\noindent
and the pair of processes
\begin{equation*}
\bigg(\bar{Y}(s), \bar{Z}(s)\bigg)\triangleq \bigg(Y_{t,x}(s\wedge\tau), 1_{[t,\tau]}(s)Z_{t,x}(s)\bigg), \quad t\leq s\leq t+\alpha,
\end{equation*}
\vskip8pt
\noindent
which can be easily seen to be  the solution of the BSDE

\begin{equation*}
\bar{Y}(s)=v\big(\tau, X_{t,x}(\tau)\big)+\int_{s\wedge\tau}^{\tau} f\Big(\theta, X_{t,x}(\theta),v\big(\theta,X_{t,x}(\theta)\big), \bar{Z}(\theta)\Big)d\theta-\int_{s}^{t+\alpha} \bar{Z}(\theta)dW(\theta).
\end{equation*}
\vskip8pt
\noindent
Moreover, an application of IKW formula implies that the pair $(\hat{Y},\hat{Z})$ of processes  given by

\begin{align*}
\hat{Y}(s)&\triangleq \phi\big(s\wedge\tau, X_{t,x}(s\wedge\tau)\big),  \\[13pt]
\hat{Z}(s)& \triangleq  1_{[t,\tau]}(s)\Big[\left(\nabla_x\phi\right)^{\top}\big(s,X_{t,x}(s)\big)\sigma\Big(s,X_{t,x}(s),v\big(s,X_{t,x}(s)\big)\Big)+\Psi^{\phi}\big(s,X_{t,x}(s)\big)\Big],
\end{align*}
\vskip8pt
\noindent
where $t\leq s\leq t+\alpha,$ is the solution of the BSDE
\begin{align*}
\hat{Y}(s)=\phi\big(\tau,X_{t,x}(\tau)\big)&-\int_{s\wedge\tau}^{\tau}\bigg[A^{\phi}\big(\theta,X_{t,x}(\theta)\big)
%&\qquad\qquad
+(L\phi)\Big(\theta, X_{t,x}(\theta),v\big(\theta,X_{t,x}(\theta)\big),\bar{Z}(\theta)\Big)\\[15pt]
&\qquad\qquad+tr\left\{\left(\nabla_x\Psi^{\phi}\right)^{\top}\big(\theta,X_{t,x}(\theta)\big)\,\sigma\Big(\theta,X_{t,x}(\theta),v\big(\theta,X_{t,x}(\theta)\big)\Big)\right\}\bigg]d\theta\\[15pt]
& -\int_{s}^{t+\alpha} \hat{Z}(\theta)dW(\theta).
\end{align*}
\vskip8pt
Setting
\begin{align*}
\hat{\beta}(s)\triangleq -\bigg[&A^{\phi}\big(s,X_{t,x}(s)\big)+(L\phi)\Big(s, X_{t,x}(s),v\big(s,X_{t,x}(s)\big),\hat{Z}(s)\Big)\\[15pt]
&+f\Big(s, X_{t,x}(s),v\big(s,X_{t,x}(s)\big),\hat{Z}(s)\Big)\\[18pt]
&+tr\left\{\left(\nabla_x\Psi^{\phi}\right)^{\top}\big(s,X_{t,x}(s)\big)\,\sigma\Big(s,X_{t,x}(s),v\big(s,X_{t,x}(s)\big)\Big)\right\}\bigg]
\end{align*}
and

\begin{align*}
\bar{\beta}(s)\triangleq -\bigg[&A^{\phi}\big(s,X_{t,x}(s)\big)+(L\phi)\Big(s,X_{t,x}(s),v\big(s,X_{t,x}(s)\big),\bar{Z}(s)\Big)\\[15pt]
&+f\Big(s,X_{t,x}(s),v\big(s,X_{t,x}(s)\big),\bar{Z}(s)\Big)\\[15pt]
&+tr\left\{\left(\nabla_x\Psi^{\phi}\right)^{\top}\big(s,X_{t,x}(s)\big)\,\sigma\Big(s,X_{t,x}(s),v\big(s,X_{t,x}(s)\big)\Big)\right\}\bigg],
\end{align*}
\vskip8pt
\noindent
we observe that
\begin{equation*}
|\bar{\beta}(s)-\hat{\beta}(s)|\leq c ||\bar{Z}(s)-\hat{Z}(s)||
\end{equation*}
\vskip5pt
\noindent
holds for some $c>0$ and $ t\leq s\leq t+\alpha.$ Therefore there exists a bounded $\mathbb{F}$-adapted  $\mathbb{R}^{1\times d}$-valued process $\gamma(\cdot)$  such that
\begin{equation*}
|\bar{\beta}(s)-\hat{\beta}(s)|= \langle\gamma(s), \bar{Z}(s)-\hat{Z}(s)\rangle, \qquad  t\leq s\leq t+\alpha.
\end{equation*}
Defining now
$$\big(\tilde{Y}(s),\tilde{Z}(s)\big)\triangleq\big(\hat{Y}(s)-\bar{Y}(s), \hat{Z}(s)-\bar{Z}(s)\big), \qquad  t\leq s\leq t+\alpha, $$
 we get that
\begin{align*}
\tilde{Y}(s)=\phi\big(\tau,X_{t,x}(\tau)\big)&-v\big(\tau,X_{t,x}(\tau)\big)+\int_{s\wedge \tau}^{\tau} \big[\hat{\beta}(\theta)+\langle\gamma(\theta), \tilde{Z}(\theta)\rangle\big]d\theta-\int_{s\wedge \tau}^{\tau} \tilde{Z}(\theta)dW(\theta).
\end{align*}
\vskip8pt
\noindent
From the proof of Theorem 1.6 in Pardoux \cite{Pardoux-96}
 we obtain that
\begin{align*}
\tilde{Y}(t)=E\bigg[\Gamma_{t,\tau} \tilde{Y}(\tau)+\int_{t}^{\tau}\Gamma_{t,s}\,\hat{\beta}(s)ds\bigg],
\end{align*}
where

\begin{align*}
\Gamma_{t,s}=\exp\left(\int_{t}^{s}\gamma(\theta)\,dW(\theta)-\frac{1}{2}\int_{t}^{s}|\gamma(\theta)|^2d\theta\right).
\end{align*}
On the other hand, the relationship  of (\ref{eqn:contrad}) and the choice  of the stopping time $\tau$ yield the inequalities:
\begin{align*}
\tilde{Y}(\tau)&\geq 0,\\[10pt]
 \hat{\beta}(\theta)&>0, \quad \forall \ \theta\in [t,\tau], \quad \tau>t,
\end{align*}
which conclude to $\tilde{Y}(t)>0.$ This means that $v(t,x)<\phi(t,x)$ which is a contradiction to our initial assumptions. \hfill$\square$

Theorem \ref{thm:class-visc} points out the need to investigate the solvability, in terms of a unique either classical or stochastic viscosity solution, of the BSPDE (\ref{equ:BSPDEu}) as an autonomous mathematical object, i.e., independently of the related system of FBSDEs (\ref{eqn:FBSDEsanss-FINITE}). However, the study of this interesting problem is beyond the scope of the present thesis and will be reported in future research.
\vskip10pt

\subsection{Stationary solutions and the infinite horizon case}

Given the unique solution of the FBSDEs system (\ref{eqn:FBSDEsanss-FINITE}) by Theorem \ref{thm:exuniqgeneral}, Proposition \ref{thm:svs} provides a stochastic viscosity solution of the associated parabolic BSPDE (\ref{equ:BSPDEu})
 for any finite time-horizon. The main task of this subsection is to explore whether this solution preserves its probabilistic properties subject to an infinite time-horizon. In particular, in case the coefficients
 of the FBSDEs (\ref{eqn:FBSDEsanss-FINITE}) were deterministic and independent of time, then  BSPDE (\ref{equ:BSPDEu}) would reduce to an elliptic PDE whose solution is the well known representing function.
%and the probability distribution of its solution would be unrelated to any chosen initial condition.
Nevertheless, since we are now considering random coefficients for the system of FBSDEs, the best to ask would be to encounter the same random field solution of the BSPDE, as regards its probability distribution, for every chosen starting point. Therefore, in what follows we shall investigate the conditions under which the stochastic viscosity solution of BSPDE (\ref{equ:BSPDEu}) becomes \emph{stationary}.

 Consider  the family of shift operators $\{\vartheta_{t}\}_{t \in {\mathbb R}_{+}}$, such that
each mapping $\vartheta_t:(\Omega, \mathcal{F})\rightarrow(\Omega, \mathcal{F}),\ t\geq 0,$  is measurable and acts as a shift on the path of the Wiener process; that is,
$$\vartheta_t\circ W(s)=W(s+t)-W(t), \quad \forall \ t,s\geq 0.$$
Moreover  assume that the random fields $b,\sigma$ and $f,$ which are the data  of the FBSDEs system (\ref{eqn:FBSDEsanss}), enjoy stationarity properties.

\begin{hypothesis}\label{hyp:stat} The random fields $b,f,\sigma$ are stationary, i.e.,
for every $t,s\geq 0$ and $(X,Y,Z) \in  \mathbb{R}^n\times \mathbb{R}\times \mathbb{R}^{1\times d}$  we have that
$\vartheta_t\circ {\mathfrak S}(s,X ,Y,Z)={\mathfrak S}(s+t, X,Y,Z)$ for
${\mathfrak S}=b,f,\sigma$.
%\begin{comment}
%\begin{align*}
%\theta_t\circ b(s,\, \cdot\, ,\,\cdot\, ,\,\cdot\,)&=b(s+t,\, \cdot\, ,\,\cdot\, ,\,\cdot\,),\\[8pt]
%\theta_t\circ \sigma(s,\, \cdot\, ,\,\cdot\,)&=\sigma(s+t,\, \cdot\, ,\,\cdot\,),\\[8pt]
%\theta_t\circ f(s,\, \cdot\, ,\,\cdot\, ,\,\cdot\,)&=f(s+t,\, \cdot\, ,\,\cdot\, ,\,\cdot\,).
%\end{align*}
%\end{comment}
\end{hypothesis}
Then we are ready to state the following result.
\begin{theorem}\label{thm:svsstat}
Suppose that the random fields $b,$ $\sigma$ and $f$ are globally continuous and satisfy Hypotheses \ref{Hyp} and \ref{hyp:stat} with $k_5=0$ and $\psi(y)=y, \ y\geq 0.$ Then the random field
$v$ of (\ref{eqn:vy}), in terms of the unique adapted solution of the FBSDEs (\ref{eqn:FBSDEsanss}), has an indistinguishable version $\bar{v}$
such that, if it admits a semimartingale
decomposition, the random field pair $(\bar{v},\Psi^{\bar{v}})$
is a stationary stochastic viscosity solution of the BSPDE (\ref{equ:BSPDEu}), i.e.,
\begin{equation}\label{eqn:stationbarv}
\vartheta_{\alpha}\circ \bar{v}(t,x)= \bar{v}(t+\alpha,x),\quad \forall \ \, (t,x)\in\mathbb{R}_{+}\times\mathbb{R}^n \quad \text{and} \quad \alpha\geq 0.
\end{equation}
\end{theorem}
\noindent
\textbf{Proof.}
For every initial conditions $(t,x)\in\mathbb{R}_{+}\times\mathbb{R}^n$ and $s\geq t,$
the inequalities in (H3) and (H4) of Hypothesis \ref{Hyp} imply that for any $(X,Y,Z) \in  \mathbb{R}^n\times \mathbb{R}\times \mathbb{R}^{1\times d}$,
\begin{align*}
 \big|b\big(s, X,Y, Z\big)\big|^2
&\leq \ 2\big|b\big(s,X, 0, 0\big)\big|^2+2\big|b\big(s, X,Y, Z\big)-b\big(s, X, 0, 0 \big)\big|^2\\[10pt]
&\leq \ 4|b(s, 0, 0, 0)|^2+4k^2+4k^2|X|^2+4k_1^2|Y|^2+4k_2^2|Z|^2,
\end{align*}
so setting $(X,Y,Z)=(X_{t,x}(s),Y_{t,x}(s),Z_{t,x}(s))$ yields the estimate
\begin{align*}
 \big|b\big(s, X_{t,x}(s),Y_{t,x}(s), Z_{t,x}(s)\big)\big|^2 \leq \ 4|b(s, 0, 0, 0)|^2&+4k^2+4k^2|X_{t,x}(s)|^2\\[10pt]
&+4k_1^2|Y_{t,x}(s)|^2+4k_2^2||Z_{t,x}(s)||^2.
\end{align*}
Similarly for  any $(X,Y,Z) \in  \mathbb{R}^n\times \mathbb{R}\times \mathbb{R}^{1\times d}$,
\begin{align*}
 \big|f\big(s, X,Y,\big)\big|^2
&\leq  \ 2\big|f\big(s, 0,Y, 0\big)\big|^2+2\big|f\big(s, X,Y, Z\big)-f\big(s, 0, Y, 0\big)\big|^2\\[10pt]
&\leq  \ 4|f(s, 0, 0, 0)|^2+4|Y|^2+4c_1^2|X|^2+4c_2^2|Z|^2,
\end{align*}
so setting $(X,Y,Z)=(X_{t,x}(s),Y_{t,x}(s),Z_{t,x}(s))$ yields the estimate
\begin{equation*}
 \big|f\big(s, X_{t,x}(s),Y_{t,x}(s), Z_{t,x}(s)\big)\big|^2
\leq  \ 4|f(s, 0, 0, 0)|^2+4|Y_{t,x}(s)|^2+4c_1^2|X_{t,x}(s)|^2+4c_2^2||Z_{t,x}(s)||^2.
\end{equation*}
Finally,  we have that  for  any $(X,Y,Z) \in  \mathbb{R}^n\times \mathbb{R}\times \mathbb{R}^{1\times d}$,
\begin{align*}
 \big|\big|\sigma\big(s, X,Y, Z\big)\big|\big|^2 &\leq   \ 2\big|\big|\sigma\big(s, 0, 0, 0\big)\big|\big|^2+2\big|\big|\sigma\big(s, X,Y, Z\big)-\sigma\big(s, 0, 0, 0 \big)\big|\big|^2\\[10pt]
&\leq \ 2||\sigma(s, 0, 0, 0)||^2+2k_3|X|^2+2k_4|Y|^2+2k_5|Z|^2,
\end{align*}
so that
\begin{align*}
 \big|\big|\sigma\big(s, X_{t,x}(s),Y_{t,x}(s), Z_{t,x}(s)\big)\big|\big|^2 \leq   2||\sigma(s, 0, 0, 0)||^2&+2k_3|X_{t,x}(s)|^2\\[8pt]
 &+2k_4|Y_{t,x}(s)|^2+2k_5|Z_{t,x}(s)|^2.
\end{align*}
Due to the fact that $(X_{t,x},Y_{t,x}, Z_{t,x})$ is the unique solution of the infinite time-horizon FBSDEs system (\ref{eqn:FBSDEsanss}) and given
the condition (H5) of Hypothesis \ref{Hyp}, we have that the process triplet
$\Big(b\big(\cdot, X_{t,x}(\cdot),Y_{t,x}(\cdot), Z_{t,x}(\cdot)\big), f\big(\cdot, X_{t,x}(\cdot),Y_{t,x}(\cdot), Z_{t,x}(\cdot)\big), \sigma\big(\cdot, X_{t,x}(\cdot),Y_{t,x}(\cdot), Z_{t,x}(\cdot)\big) \Big)$
is also of class $\mathbb{M}_{\lambda}^2(t,\infty;\mathbb{R}^n\times\mathbb{R}\times\mathbb{R}^{1\times d}).$

Instead of (\ref{eqn:FBSDEsanss}), we consider  its  finite horizon approximation  in $[t,T]:$
\begin{align}\label{eqn:FBSDEsanssequiv}
%\begin{split}\label{eqn:FBSDEsanssequiv}
X_{t,x}(T)&=x+\int_{t}^{T}b\big(s,X_{t,x}(s),Y_{t,x}(s), Z_{t,x}(s)\big)ds+\int_{t}^{T}\sigma\big(s,X_{t,x}(s), Y_{t,x}(s)\big)dW(s),\nonumber\\ \\
Y_{t,x}(t)&=Y_{t,x}(T)+\int_{t}^{T}f\big(s,X_{t,x}(s), Y_{t,x}(s), Z_{t,x}(s)\big)ds-\int_{t}^{T} Z_{t,x}(s)dW(s),\nonumber
%\end{split}
\end{align}
\vskip5pt
\noindent
where $\lim_{T\rightarrow\infty}e^{\lambda T} Y_{t,x}^2(T)=0$,  so that $Y_{t,x}(\cdot)$ may satisfy the infinite horizon problem (\ref{eqn:FBSDEsanss}).  We will vary $(t,x,T)$ and consider \eqref{eqn:FBSDEsanssequiv} as a family of finite horizon FBSDEs in the following sense. For the specific choice of $(t,x,T)$, the FBSDEs system of \eqref{eqn:FBSDEsanssequiv}  has as solution a triplet of processes $\left(X_{t,x},Y_{t,x},Z_{t,x}\right)$ that can be interpreted as an approximation to the solution of the infinite horizon FBSDEs  (\ref{eqn:FBSDEsanss}), with the forward component $X_{t,x}(\cdot)$ starting at $t$ at position $x$ and running up to final time $T$. Then Hypothesis \ref{hyp:stat} and an application of the operator $\vartheta_{\alpha}$ for any $\alpha\geq 0$
to the FBSDEs  (\ref{eqn:FBSDEsanssequiv}) yield the system
\begin{align}\label{eqn:FBSDEsanssequivtheta}
\theta_{\alpha}\circ X_{t,x}(T)&=x+\int_{t+\alpha}^{T+\alpha}b\big(s,\theta_{\alpha}\circ X_{t,x}(s-\alpha), \theta_{\alpha}\circ Y_{t,x}(s-\alpha), \theta_{\alpha}\circ Z_{t,x}(s-\alpha)\big)ds\nonumber\\[10pt]
&\qquad +\int_{t+\alpha}^{T+\alpha}\sigma\big(s,\theta_{\alpha}\circ X_{t,x}(s-\alpha), \theta_{\alpha}\circ Y_{t,x}(s-\alpha)\big)dW(s),\nonumber\\[10pt]
\theta_{\alpha}\circ Y_{t,x}(t)&=\theta_{\alpha}\circ Y_{t,x}(T)+\int_{t+\alpha}^{T+\alpha}f\big(s,\theta_{\alpha}\circ X_{t,x}(s-\alpha), \theta_{\alpha}\circ Y_{t,x}(s-\alpha), \theta_{\alpha}\circ Z_{t,x}(s-\alpha)\big)ds\nonumber\\[10pt]
&\qquad\qquad\qquad-\int_{t+\alpha}^{T+\alpha} \theta_{\alpha}\circ Z_{t,x}(s-\alpha)dW(s),  \\[10pt]
\lim_{T\rightarrow\infty}&e^{\lambda (T+\alpha)}(\theta_{\alpha}\circ Y_{t,x}^2(T))=0.\nonumber
\end{align}
%\vskip5pt
\noindent
We now consider the family of finite horizon FBSDEs \eqref{eqn:FBSDEsanssequiv} for the choice $(t+\alpha,x,T+\alpha)$ to obtain the triplet of processes
$(X_{t+\alpha,x},Y_{t+\alpha,x}, Z_{t+\alpha,x}),$ which is an approximation to the solution of the infinite horizon FBSDEs (\ref{eqn:FBSDEsanss}) if we start at $t+\alpha,$ setting $ X_{t,x}(\cdot)$ equal to $x,$ and  run the system for horizon $T+\alpha$.  For this choice   (\ref{eqn:FBSDEsanssequiv}) can be written as
\begin{align}\label{eqn:FBSDEsanssequivxi}
X_{t+\alpha,x}(T+\alpha)&=x+\int_{t+\alpha}^{T+\alpha}b\big(s,X_{t+\alpha,x}(s),Y_{t+\alpha,x}(s), Z_{t+\alpha,x}(s)\big)ds\nonumber\\[10pt]
&\qquad+\int_{t+\alpha}^{T+\alpha}\sigma\big(s,X_{t+\alpha,x}(s), Y_{t+\alpha,x}(s)\big)dW(s),\\[12pt]
Y_{t+\alpha,x}(t+\alpha)&=Y_{t+\alpha,x}(T+\alpha)+\int_{t+\alpha}^{T+\alpha}f\big(s,X_{t+\alpha,x}(s), Y_{t+\alpha,x}(s), Z_{t+\alpha,x}(s)\big)ds\nonumber\\[12pt]
&\qquad\qquad\qquad-\int_{t+\alpha}^{T+\alpha} Z_{t+\alpha,x}(s)dW(s),  \nonumber\\[15pt]
\lim_{T\rightarrow\infty}&e^{\lambda (T+\alpha)} Y_{t+\alpha,x}^2(T+\alpha)=0.\nonumber
\end{align}

Therefore,  a comparison of the FBSDEs systems  (\ref{eqn:FBSDEsanssequivtheta}) and (\ref{eqn:FBSDEsanssequivxi}),
the uniqueness of solution to the  FBSDEs system (\ref{eqn:FBSDEsanss}) combined with its continuity according to the initial parameters
(cf. Corollary \ref{col:cdsp}), and a perfection procedure (cf. Arnold \cite{Arnold}
and Arnold and Scheutzow \cite{Arnold-Scheutzow})
imply that
 \begin{align*}
\vartheta_{\alpha}\circ X_{t,x}(s)&=X_{t+\alpha,x}(s+\alpha), \quad \vartheta_{\alpha}\circ Y_{t,x}(s)=Y_{t+\alpha,x}(s+\alpha)\quad \text{and}\\[10pt]
\vartheta_{\alpha}\circ Z_{t,x}(s)&=Z_{t+\alpha,x}(s+\alpha), \quad \forall \ \ \alpha\geq 0, \quad s\geq t, \quad x\in\mathbb{R}^{n};
\end{align*}
in other words, the process triplet $(X_{t,x},Y_{t,x}, Z_{t,x})$ is a stationary solution of the FBSDEs system (\ref{eqn:FBSDEsanss}).
As a consequence from (\ref{eqn:vy}),  for any  $(t,x)\in\mathbb{R}_{+}\times\mathbb{R}^n$  the random field $v$ satisfies that
$$\vartheta_{\alpha}\circ v(t,x)=v(t+\alpha,x), \quad \forall \ \ \alpha\geq 0. $$
Due though to Proposition \ref{thm:svs} this random field belongs to the class of  $C_{\mathbb{F}}\big([0,T];$ $\mathbb{L}^2(\Omega;C(\mathbb{R}^{n};\mathbb{R}))\big)$,
thus we can get an indistinguishable version of $v,$ denoted by $\bar{v}$, which satisfies the stationarity relationship of (\ref{eqn:stationbarv})
and the remaining assertion of the theorem.
\hfill$\square$

\begin{remark}\label{rem:compar-viscos}
One should not fail to notice that Theorem \ref{thm:comparison} implies a comparison property for the stochastic viscosity solutions
of the BSPDE (\ref{equ:BSPDEu}). In particular, let the random fields $b_i,$ $f_i,$ $i=1,2,$ and $\sigma$ be globally continuous, satisfy Hypothesis \ref{Hyp}
with $k_5=0,$ as well as the relationships of (\ref{eqn:bifisigma}), and consider the BSPDEs

\begin{align}\label{equ:BSPDEui}
v_i(t,x)=v_i(T,x)&+\int_{t}^{T}\bigg[(L_iv_i)\Big(s,x,v_i(s,x),\left(\nabla_x v_i\right)^{\top}(s,x)\,\sigma\big(s,x,v_i(s,x)\big)+\Psi^{v_i}(s,x)\Big)\nonumber\\[12pt]
& \qquad \ +f_i\Big(s,x,v_i(s,x),\left(\nabla_x v_i\right)^{\top}(s,x)\,\sigma\big(s,x,v_i(s,x)\big)+\Psi^{v_i}(s,x)\Big)\nonumber\\[12pt]
& \qquad \ +tr\left\{\left(\nabla_x \Psi^{v_i}\right)^{\top}(s,x)\,\sigma\big(s,x,v_i(s,x)\big)\right\}\bigg]ds\\[12pt]
&-\int_{t}^{T}\Psi^{v_i}(s,x)dW(s), \quad \forall\quad (t,x)\in[0,T]\times\mathbb{R}^{n},\nonumber
\end{align}
 where
\begin{equation*}
(L_iv_i)(t,\bar{x},y,z)\triangleq \frac{1}{2}\sum_{k,j=1}^{n} \big((\sigma\sigma^{\top})(t,\bar{x},y)\big)_{k,j}\, \frac{\partial^2v_i}{\partial x_k\partial x_j} (t,\bar{x})+\big\langle b_i(t,\bar{x},y,z),\nabla_x v_i(t,\bar{x})\big\rangle
\end{equation*}
for all $(t,\bar{x},y,z)\in[0,T)\times\mathbb{R}^{n}\times\mathbb{R}\times\mathbb{R}^{1\times d}$ and $i=1,2.$
If $(X_{i,t,x}, Y_{i,t,x}, Z_{i,t,x})$
is the unique adapted solution of the infinite horizon FBSDEs
\begin{align}\label{eqn:FBSDEsanssi}
dX_{i,t,x}(s)&=b_i\big(s,X_{i,t,x}(s),Y_{i,t,x}(s), Z_{i,t,x}(s)\big)ds+\sigma\big(s,X_{i,t,x}(s), Y_{i,t,x}(s)\big)dW(s),\ s>t,\nonumber\\[10pt]
dY_{i,t,x}(s)&=-f_i\big(s,X_{i,t,x}(s), Y_{i,t,x}(s), Z_{i,t,x}(s)\big)ds+ Z_{i,t,x}(s)dW(s), \quad s\geq t,  \\[10pt]
X_{i,t,x}(t)&=x,\nonumber
\end{align}
and the random field
\begin{equation*}
 v_i(t,x)\triangleq Y_{i,t,x}(t),\qquad \forall\ (t,x)\in[0,T]\times \mathbb{R}^n,
 \end{equation*}
 obtains a semimartingale decomposition, then Proposition \ref{thm:svs} yields that the latter is of class $C_{\mathbb{F}}\big([0,T];\mathbb{L}^2(\Omega;C(\mathbb{R}^{n};\mathbb{R}))\big)$  and
 the random field pair $(v_i,\Psi^{v_i})$ is a stochastic viscosity solution of BSPDE (\ref{equ:BSPDEui}) for $i=1,2$.
By Theorem \ref{thm:comparison} we obtain directly that $ v_1(t,x)\geq  v_2(t,x)$ for every  $(t,x)\in[0,T]\times \mathbb{R}^n.$
Clearly, by a similar analysis subject now to Theorem \ref{thm:svsstat}, this comparison property holds also for the stationary stochastic viscosity solutions of BSPDE (\ref{equ:BSPDEu}).
\end{remark}

In the higher-dimensional case of $m>1,$ it is straightforward to see that the reasoning deployed in this section still applies to every coordinate of the
resulting $\mathbb{R}^m$-valued representing random field $v$ of (\ref{FSSC}).
In particular, the BSPDE of (\ref{equ:BSPDEu}) is replaced by the system of quasilinear BSPDEs of the form
\begin{align*} %\label{eqn:BSPDEvecu}
v^k(t,x)=v^k(T,x)&+\int_{t}^{T}\bigg[(Lv^k)\Big(s,x,v(s,x),\left(\nabla_x v\right)^{\top}(s,x)\,\sigma\big(s,x,v(s,x)\big)+\Psi^v(s,x)\Big)\nonumber\\[12pt]
& \qquad \quad  \ +f\Big(s,x,v(s,x),\left(\nabla_x v\right)^{\top}(s,x)\,\sigma\big(s,x,v(s,x)\big)+\Psi^{v}(s,x)\Big)\nonumber
\end{align*}
\begin{align*}
& \qquad \qquad \qquad \qquad +tr\left\{\left(\nabla_x \Psi^{v^k}\right)^{\top}(s,x)\,\sigma\big(s,x,v(s,x)\big)\right\}\bigg]ds\\[12pt]
&\qquad\qquad+\int_{t}^{T}\Psi^{v^k}(s,x)dW(s), \quad \forall \ \, k=1,...,m,  \ \, (t,x)\in[0,T]\times\mathbb{R}^{n}.\nonumber
\end{align*}

\vskip25pt
\section{Infinite Horizon Stochastic Maximum Principle Of Optimal Control To Random Forward - Backward Systems} \label{sec:ihmp}

In Section \ref{sec:ssa} we derived existence and uniqueness results for an adapted solution of the infinite horizon fully coupled system of random FBSDEs (\ref{GFBSDEs}),
cf. Theorem \ref{thm:exuniqgeneral} and Remark \ref{rem:exisuniqsigma}, whose study was motivated by several future expectations models of economic theory in Section \ref{sec:rem}.
In the present section we shall formulate these models to motivate now the investigation of the associated
stochastic control problem of such FBSDEs systems. Then,  moving along the lines of
 Maslowski and Veverka \cite{Veverka},
 we shall establish a sufficient \emph{stochastic maximum principle}
for the \emph{exponentially} $\lambda$\emph{-weighted} control problem that arises from the multidimensional fully coupled FBSDEs (\ref{GFBSDEs}) with a general control domain.
Making in particular use of the duality methodology,
we shall define the \emph{generalized Hamiltonian} which will allow us to state the \emph{adjoint system} of FBSDEs. 
\vskip15pt
\subsection{Stochastic optimal control of future expectations models}\label{sub:socfem}
In Subsections \ref{KRUGMAN-SECTION}-\ref{sub:consolrate} we saw that if we extend the dynamics of several well known rational expectations models
to accommodate random coefficients, so as to reflect better the course of real economy, then they obey the general stochastic saddlepoint system of (\ref{REM}),
which due to Proposition \ref{Prop:FBSDEsToREM} is equivalent to the infinite horizon random FBSDEs system (\ref{FBSDEREM}).
However, the monetary authorities have the right to intervene in these dynamics such as to maintain the corresponding economic variables either
 inside a predetermined zone or as close as possible
to a given target, but whenever they do so they also suffer intervention costs. This intervention bevavior can be quantified through a control
variable $u(\cdot)$ that is incorporated by the coefficients of the above systems;
in fact, the controlled saddlepoint system of (\ref{REM}) has now the form:
\begin{align}\label{CREM}
dX(t)&=b\big(t,X(t),Y(t),u(t)\big)dt+\sigma\big(t,X(t),Y(t),u(t)\big)dW(t),\quad t>0, \nonumber\\[13pt]
Y(t)&=E\left[\int_{t}^{\infty} e^{-\int_{t}^{s}r(\theta)d\theta}g\big(s,X(s),Y(s),u(s)\big)ds\Big|\mathcal{F}(t)\right],\quad t\geq 0,\\[13pt]
X(0)&=x\in\mathbb{R}^{n}, \nonumber
\end{align}
which is equivalent to the controlled infinite horizon FBSDEs

\begin{align}\label{CFBSDEREM}
dX(t)&=b\big(t,X(t),Y(t),u(t)\big)dt+\sigma\big(t,X(t),Y(t), u(t)\big)dW(t), \quad t> 0,\nonumber\\[10pt]
dY(t)&=\left[-g\big(t, X(t),Y(t), u(t)\big)+r(t)Y(t)\right]dt+Z(t)dW(t), \quad t\geq 0,\\[10pt]
X(0)&=x\in\mathbb{R}^{n}.\nonumber
\end{align}

Here we assume that the monetary authorities regulate the control process $u(\cdot)$ of the saddlepoint system (\ref{CREM}) so that its variable
pair $(X,Y)$ is kept as close as possible to a given target pair process $(c_1, c_2)$ rather than within a pre-described band. Since though
 this is a costly procedure this regulation activity must be the minimum possible.Therefore, the control process $u(\cdot)$ could be chosen
 as such as to minimize the \emph{quadratic} $\lambda$-weighted \emph{cost functional}
 \begin{align}\label{COSREM}
J_0\big(x;u(\cdot)\big)\triangleq E\Bigg[&\int_0^{\infty} e^{\lambda \, t} \Big\{ P(t)\big[X(t)-c_1(t)\big]^2+ Q(t)\big[Y(t)-c_2(t)\big]^2
+R(t)u^2(t)\Big\}dt \\[8pt]
 & +N\big[Y(0)-c_2(0)\big]^2\Bigg],\nonumber
\end{align}
where
$P(\cdot), Q(\cdot), R(\cdot)>0$ are weight processes representing the cost of deviation from the target pair and the cost of intervention,
respectively,  $N$ is a weight that corresponds to the initial value of variable $Y(\cdot),$ and  $\lambda\in\mathbb{R}$ such as the exponential term ensures
the convergence of the cost functional.

Variations of both Krugman and Dornbusch rational expectations models, presented in Subsections \ref{KRUGMAN-SECTION} and \ref{DORNBUSCH-SECTION} respectively, may be used to provide simple illustrations of the above stochastic optimal control
problem. In particular, the monetary authorities regulate the exogenously given domestic money supply $m(\cdot)$ in order to keep
the domestic price of foreign exchange  $s(\cdot)$ close to a given target, subject to an inflation cost represented by the functional of (\ref{COSREM}). This version of Krugman's model
has been treated in the finite time-horizon case by Yannacopoulos \cite{Yannacopoulos05}
as a stochastic control problem of a BSDE, since its FSDE is decoupled from the backward process. However the corresponding Dornbusch's model
is an optimization problem of a strongly coupled FBSDEs system. Several other formulations of the exchange rate control problem may be found
in Jeanblanc-Picqu\'{e} \cite{Jeanblanc},
Mundaca and {\O}ksendal \cite{Mundaca},
Cadenillas and Zapatero \cite{Cadenillas}
 and references therein.

\vskip10pt
\subsection{Stochastic maximum principle of infinite horizon random FBSDEs}
 Let $U$ be a nonempty subset of $\mathbb{R}^k$ and $\lambda \in \mathbb{R}$. In light of the optimal control problem
 we described in the previous subsection through (\ref{CREM})-(\ref{COSREM}), for any \emph{admissible control} process $u(\cdot)$ that belongs to
 the set
\[\mathcal{U}_{\lambda}\triangleq\left\{u:\mathbb{R}_{+}\times \Omega \rightarrow \mathbb{R}^k \, \Big{|}\, u(\cdot)\in \mathbb{M}_{\lambda}^2(0,\infty; U)\right\},\]
we consider the stochastic control problem whose \emph{state variables} are described by the infinite horizon FBSDEs
\begin{align}\label{FBSDEu}
dX(t)&=b\big(t,X(t),Y(t),Z(t), u(t)\big)dt+\sigma\big(t,X(t),Y(t),Z(t), u(t)\big)dW(t), \ t>0,\nonumber\\[10pt]
dY(t)&=-f\big(t,X(t),Y(t),Z(t), u(t)\big)dt+Z(t)dW(t),\quad t\geq 0,\\[10pt]
X(0)&=x\in \mathbb{R}^n \nonumber,
\end{align}
and the \emph{cost functional} is defined by the total exponentially $\lambda$-weighted average
\begin{equation}\label{eqn:cost}
J\big(x;u(\cdot)\big)\triangleq E\left[\int_0^{\infty} e^{\lambda \, t} \, h\big(t,X(t),Y(t),Z(t), u(t)\big)dt +a\big(Y(0)\big)\right],
\end{equation}
where $b: \mathbb{R}_{+}\times \mathbb{R}^{n}\times \mathbb{R}^{m}\times \mathbb{R}^{m\times d}\times \mathbb{R}^k\times \Omega \rightarrow \mathbb{R}^n , \
\sigma:\mathbb{R}_{+}\times \mathbb{R}^{n}\times \mathbb{R}^{m}\times \mathbb{R}^{m\times d}\times \mathbb{R}^k\times\Omega \rightarrow \mathbb{R}^{n\times d}, \
 f:\mathbb{R}_{+}\times \mathbb{R}^{n}\times \mathbb{R}^{m}\times \mathbb{R}^{m\times d}\times \mathbb{R}^k\times \Omega \rightarrow \mathbb{R}^m,
 \ h: \mathbb{R}_{+}\times \mathbb{R}^{n}\times \mathbb{R}^{m}\times \mathbb{R}^{m\times d}\times \mathbb{R}^k\times \Omega \rightarrow \mathbb{R} $
%$b, \, \sigma, \, f, \, g$
and
$a:\mathbb{R}^{m}\times \Omega\rightarrow \mathbb{R}$ are given random fields.
Due to matters of simplicity we use the same notation for the coefficients of both systems of FBSDEs (\ref{GFBSDEs}) and (\ref{FBSDEu}),
and due to matters of consistency we make the following assumption that is an extension of Hypothesis \ref{Hyp}.
\begin{assumption}\label{ass:control}
We impose the following conditions to the coefficients of FBSDEs (\ref{FBSDEu}).

(A1) For every $u\in\mathbb{R}^k$ the random fields $b(\theta,u),$ $\sigma(\theta,u)$ and $f(\theta,u)$
satisfy Hypothesis \ref{Hyp} with respect to $\theta=(t,x,y,z),$ where the last condition of (H5) is replaced by
\[ \Big(b\big(\cdot,0,0,0,u(\cdot)\big),f\big(\cdot,0,0,0,u(\cdot)\big),\sigma\big(\cdot,0,0,0,u(\cdot)\big)\Big)\in \mathbb{M}_{\lambda}^2\big(0,\infty; \mathbb{R}^n\times\mathbb{R}^m\times \mathbb{R}^{n\times d}\big) \]
for any $u(\cdot)\in\mathcal{U}_{\lambda}.$

(A2) The \emph{penalization} random fields $h$ and $a$ are chosen such that the cost functional of (\ref{eqn:cost}) converges for every admissible control, i.e.,
\[\big|J\big(x;u(\cdot)\big)\big|<\infty, \quad \forall \ \,  u(\cdot)\in \mathcal{U}_{\lambda}.\]
In addition, the mapping $y\mapsto a(y)$ is a convex function.

(A3) The random fields $b, \, \sigma,\, f, \, h$ and $a$ are continuously differentiable with respect to $(x,y,z)$ with bounded first-order partial derivatives.
\end{assumption}

Given any process $u(\cdot)\in\mathcal{U}_{\lambda},$ Assumption \ref{ass:control} and Theorem \ref{thm:exuniqgeneral} ensure the existence of a
unique adapted solution triplet $(X, Y, Z)$ which solves the infinite horizon FBSDEs system of (\ref{FBSDEu}).
Our goal is to minimize the cost functional of (\ref{eqn:cost}) over all admissible controls $u(\cdot)\in\mathcal{U}_{\lambda};$ that is, to characterize the
\emph{optimal control} process $u^\ast(\cdot)$ that attains the infimum of the \emph{value function}
\begin{equation}\label{vf}
V(x)\triangleq \inf_{u(\cdot)\in\mathcal{U}_{\lambda}} J\big(x;u(\cdot)\big).
\end{equation}

To proceed with our analysis, we define the \emph{Hamiltonian random field} $H: \mathbb{R}_{+}\times \mathbb{R}^{n}\times \mathbb{R}^{m}\times \mathbb{R}^{m\times d}\times \mathbb{R}^k\times \mathbb{R}^{m}\times \mathbb{R}^{n}\times \mathbb{R}^{n\times d}\times \Omega\rightarrow\mathbb{R}$ that is associated with the stochastic control problem (\ref{FBSDEu})-(\ref{vf})  as

\begin{align}\label{eqn:hamilt}
H(t,x,y,z,u,p,q,k)\triangleq & \langle b(t,x,y,z,u), q\rangle+tr\big\{\sigma^{\top}(t,x,y,z,u) k\big\}\nonumber\\[10pt]
&\, - \langle f(t,x,y,z,u),p \rangle+h(t,x,y,z,u)\\[10pt]
& +\lambda \langle x,q\rangle +\lambda \langle y,p\rangle.\nonumber
\end{align}
The last two terms play the role of the Lyapunov function, used also in Maslowski and Veverka \cite{Veverka},
needed for the employment of the duality theory to the above optimal control problem; cf. Theorem \ref{thm:maxpri} below. Under Assumption \ref{ass:control}, the Hamiltonian $H$ is differentiable with respect to $(x, y, z)$ and thanks to Theorem \ref{thm:exuniqgeneral} there exists a unique adapted solution triplet $(p, q, k)$
 of the corresponding \emph{adjoint system} of FBSDEs  given by
\begin{align}\label{FBSDEpq}
dp(t)&=-\nabla_y H\big(t, X(t), Y(t), Z(t), u(t), p(t), q(t), k(t)\big) dt \nonumber\\[10pt]
&\quad-\nabla_z H\big(t, X(t), Y(t), Z(t), u(t), p(t), q(t), k(t)\big) dW(t), \quad t>0,\\[10pt]
dq(t)&=-\nabla_x H\big(t, X(t), Y(t), Z(t), u(t), p(t), q(t), k(t)\big) dt +k(t)dW(t), \quad t\geq 0, \nonumber\\[10pt]
p(0)&=-\nabla_y a\big(Y(0)\big).\nonumber
\end{align}
The result that follows constitutes a \emph{sufficient} stochastic maximum principle of the above
control problem.
\begin{theorem}\label{thm:maxpri}
Let Assumption \ref{ass:control} hold. Let also $\hat{u}(\cdot)\in \mathcal{U}_{\lambda}$ be an admissible control process
and $(\hat{X},\hat{Y},\hat{Z})$ be the solution of the corresponding controlled infinite horizon FBSDEs system (\ref{FBSDEu}) that belongs to $ \mathbb{M}_{\lambda}^2\big(0,\infty; \mathbb{R}^n\times\mathbb{R}^m\times \mathbb{R}^{m\times d}\big)$. Furthermore, let $(\hat{p},\hat{q},\hat{k})\in  \mathbb{M}_{\lambda}^2\big(0,\infty; \mathbb{R}^m\times\mathbb{R}^n\times \mathbb{R}^{n\times d}\big)$ be the solution of the associated adjoint infinite horizon FBSDEs system  (\ref{FBSDEpq}) such that for every  $ t\geq 0$
\begin{align*}
&(i) \ H\big(t,\hat{X}(t), \hat{Y}(t), \hat{Z}(t), \hat{u}(t), \hat{p}(t), \hat{q}(t), \hat{k}(t)\big)\\[8pt]
&\qquad \qquad\qquad\qquad\qquad\quad =\min_{u(\cdot)\in \mathcal{U}_{\lambda}} H\big(t,\hat{X}(t), \hat{Y}(t), \hat{Z}(t), u(t), \hat{p}(t), \hat{q}(t), \hat{k}(t)\big),
\end{align*}
\begin{align*}
(ii) \ (x,y,z,u)\mapsto H(t,x,y,z,u, \hat{p}(t), \hat{q}(t), \hat{k}(t)\big) \ \text{is a convex mapping.\phantom{grfjgnjtu}}
\end{align*}
Then $u^{\ast}(\cdot)=\hat{u}(\cdot)$ is the optimal control process for the infinite horizon stochastic control problem (\ref{FBSDEu})-(\ref{vf}), that is
\[V(x)=J\big(x;\hat{u}(\cdot)\big).\]
\end{theorem}
\noindent
\textbf{Proof.} Consider an arbitrary admissible stochastic control process $u(\cdot)\in \mathcal{U}_{\lambda}$ and the solution triplet $(X,Y,Z)$ of the respective infinite horizon FBSDEs  (\ref{FBSDEu}) that belongs to $ \mathbb{M}_{\lambda}^2\big(0,\infty; \mathbb{R}^n\times\mathbb{R}^m\times \mathbb{R}^{m\times d}\big)$. In order to prove that $\hat{u}(\cdot)$ is the optimal control process
it suffices to show that 
$J\big(x,\hat{u}(\cdot)\big)-J\big(x,u(\cdot)\big)\leq 0.$ Henceforth, making use of the definitions (\ref{eqn:cost}) and (\ref{eqn:hamilt}),
the convergence of the cost functional $J,$ and Fubini's theorem,  we have that
\begin{align*}
&\quad \ \, J\big(x,\hat{u}(\cdot)\big)-J\big(x,u(\cdot)\big)\\[15pt]
&=E\bigg\{\int_{0}^{\infty} e^{\lambda\,t}\left[h\big(t,\hat{X}(t), \hat{Y}(t), \hat{Z}(t), \hat{u}(t)\big)-h\big(t,X(t), Y(t), Z(t), u(t)\big)\right]dt\\[10pt]
&\qquad\quad+a\big(\hat{Y}(0)\big)-a\big(Y(0)\big)\bigg\}\\[15pt]
&=E\bigg\{\int_{0}^{\infty}  e^{\lambda\,t} \Big[H\big(t,\hat{X}(t), \hat{Y}(t), \hat{Z}(t), \hat{u}(t), \hat{p}(t), \hat{q}(t), \hat{k}(t) \big)\\[15pt]
&\qquad\qquad\quad-H\big(t,X(t), Y(t), Z(t), u(t),\hat{p}(t), \hat{q}(t), \hat{k}(t)\big)\\[15pt]
&\qquad \quad +\big\langle b\big(t,X(t), Y(t), Z(t), u(t)\big)-b\big(t,\hat{X}(t), \hat{Y}(t), \hat{Z}(t), \hat{u}(t)\big),\hat{q}(t)\big\rangle\\[15pt]
&\qquad \ \ \, +tr\Big\{\Big[\sigma^\top\big(t,X(t), Y(t), Z(t), u(t)\big)-\sigma^\top\big(t,\hat{X}(t), \hat{Y}(t), \hat{Z}(t), \hat{u}(t)\big)\Big]\hat{k}(t)\Big\}\\[15pt]
&\qquad \quad  +\big \langle f\big(t,\hat{X}(t), \hat{Y}(t), \hat{Z}(t), \hat{u}(t)\big)- f\big(t,X(t), Y(t), Z(t), u(t)\big),\hat{p}(t) \big\rangle\\[15pt]
&\qquad \quad +\lambda \big\langle X(t)-\hat{X}(t),\hat{q}(t)\big\rangle +\lambda \big\langle Y(t)-\hat{Y}(t),\hat{p}(t)\big\rangle\Big]dt
%\\[8pt]
%&\qquad \quad
+a\big(\hat{Y}(0)\big)-a\big(Y(0)\big)\bigg\}\\[10pt]
&=E\left[\int_{0}^{\infty} I(t)dt +a\big(\hat{Y}(0)\big)-a\big(Y(0)\big)\right]\\[15pt]
&=\lim_{T\rightarrow\infty}E\left[\int_{0}^{T} I(t)dt+a\big(\hat{Y}(0)\big)-a\big(Y(0)\big)\right],
\end{align*}
\vskip5pt
\noindent
where $I(t)$ is the integrant of the Lebesque integral preceding  the next-to-last equality.

On the other hand, by employing the Cauchy-Schwarz inequality together with an elementary algebraic inequality, we obtain that
\begin{align*}
\Big|E\left [e^{\lambda\, T}\big \langle \hat{X}(T)-X(T), \hat{q}(T)\big\rangle\right]\Big|\leq \frac{1}{2}E\left[e^{\lambda\, T} |\hat{X}(T)-X(T)|^2\right]+ \frac{1}{2}E\left[e^{\lambda\, T} |\hat{q}(T)|^2\right]
\rightarrow 0
\end{align*}
\vskip5pt
\noindent
as $T\rightarrow \infty$ since $X(\cdot), \hat{X}(\cdot)$ and $\hat{q}(\cdot)$ are of class $ \mathbb{M}_{\lambda}^2\big(0,\infty; \mathbb{R}^n\big).$ Furthermore, in view of the left-hand-side
of the above inequality and taking into account that expectations of square integrable stochastic integrals vanish, It\^{o}'s formula on $[0,T]$ implies that

\begin{align*}
E&\left[e^{\lambda\, T}\langle\hat{X}(T)-X(T),\hat{q}(T)\rangle\right]=E\int_{0}^{T}e^{\lambda\, t}\bigg[ \lambda \langle\hat{X}(t)-X(t),\hat{q}(t)\rangle\\[20pt]
&+ \big\langle b\big(t,\hat{X}(t), \hat{Y}(t), \hat{Z}(t), \hat{u}(t)\big)-b\big(t,X(t), Y(t), Z(t), u(t)\big),\hat{q}(t)\big\rangle \\[20pt]
&+tr\left\{\Big[\sigma^\top\big(t,\hat{X}(t), \hat{Y}(t), \hat{Z}(t), \hat{u}(t)\big)-\sigma^\top\big(t,X(t), Y(t), Z(t), u(t)\big)\Big]\hat{k}(t)\right\}\\[20pt]
&-\big\langle\hat{X}(t)-X(t),\nabla_x H\big(t,\hat{X}(t), \hat{Y}(t), \hat{Z}(t), \hat{u}(t), \hat{p}(t), \hat{q}(t), \hat{k}(t) \big)\big\rangle\bigg]dt.
\end{align*}
Following a similar reasoning, we may deduce that
\begin{align*}
E\left [e^{\lambda\, T}\big \langle \hat{Y}(T)-Y(T), \hat{p}(T)\big\rangle\right]
  \overset{T\rightarrow \infty}{\longrightarrow}0
\end{align*}
and
\begin{align*}
& \ E\left[e^{\lambda\, T}\langle \hat{Y}(T)-Y(T),\hat{p}(T)\rangle\right]\\[16pt]
=& \ E\bigg\{-\big\langle \hat{Y}(0)-Y(0),\nabla_y a\big(\hat{Y}(0)\big)\big\rangle\\[20pt]
&\ \ \ \quad+\int_{0}^{T}e^{\lambda\, t}\bigg[\lambda \langle \hat{Y}(t)-Y(t),\hat{p}(t)\rangle\\[17pt]
& \ \ \ \quad\qquad\qquad\quad+ \big\langle f\big(t,X(t), Y(t), Z(t), u(t)\big)-f\big(t,\hat{X}(t), \hat{Y}(t), \hat{Z}(t), \hat{u}(t)\big),\hat{p}(t)\big\rangle \\[20pt]
&\ \ \ \quad \qquad\qquad\quad-tr\Big\{ \big[\hat{Z}(t)-Z(t)\big]\nabla_z H\big(t,\hat{X}(t), \hat{Y}(t), \hat{Z}(t), \hat{u}(t), \hat{p}(t), \hat{q}(t), \hat{k}(t) \big)\Big\}\\[20pt]
&\ \ \ \quad \qquad\qquad\quad-\big\langle  \hat{Y}(t)-Y(t), \nabla_y H\big(t,\hat{X}(t), \hat{Y}(t), \hat{Z}(t), \hat{u}(t), \hat{p}(t), \hat{q}(t), \hat{k}(t) \big)\big\rangle\bigg]dt\bigg\}
\end{align*}
hold as well.

Finally, from all the above relationships we compute directly that

\begin{align*}
&\quad \ \, J\big(x,\hat{u}(\cdot)\big)-J\big(x,u(\cdot)\big)\\[17pt]
&=\lim_{T\rightarrow\infty}E\bigg[\int_{0}^{T} I(t)dt+e^{\lambda\, T}\langle\hat{X}(T)-X(T),\hat{q}(T)\rangle+e^{\lambda\, T}\langle \hat{Y}(T)-Y(T),\hat{p}(T)\rangle\\[17pt]
&\qquad\qquad\quad+a\big(\hat{Y}(0)\big)-a\big(Y(0)\big)\bigg]\\[17pt]
&=\lim_{T\rightarrow\infty}E\bigg\{\int_{0}^{T}  e^{\lambda\,t} \Big[H\big(t,\hat{X}(t), \hat{Y}(t), \hat{Z}(t), \hat{u}(t), \hat{p}(t), \hat{q}(t), \hat{k}(t) \big)\\[14pt]
&\qquad\qquad\quad-H\big(t,X(t), Y(t), Z(t), u(t),\hat{p}(t), \hat{q}(t), \hat{k}(t)\big)\\[10pt]
&\qquad\qquad\quad-\big\langle\hat{X}(t)-X(t),\nabla_x H\big(t,\hat{X}(t), \hat{Y}(t), \hat{Z}(t), \hat{u}(t), \hat{p}(t), \hat{q}(t), \hat{k}(t) \big)\big\rangle\\[14pt]
&\qquad\qquad\quad-\big\langle  \hat{Y}(t)-Y(t), \nabla_y H\big(t,\hat{X}(t), \hat{Y}(t), \hat{Z}(t), \hat{u}(t), \hat{p}(t), \hat{q}(t), \hat{k}(t) \big)\big\rangle\\[14pt]
&\qquad\qquad\quad-tr\Big\{ \big[\hat{Z}(t)-Z(t)\big]\nabla_z H\big(t,\hat{X}(t), \hat{Y}(t), \hat{Z}(t), \hat{u}(t), \hat{p}(t), \hat{q}(t), \hat{k}(t) \big)\Big\}\Big]dt\\[14pt]
&\qquad\qquad \ +a\big(\hat{Y}(0)\big)-a\big(Y(0)\big)-\big\langle\hat{Y}(0)- Y(0),\nabla_y a\big(\hat{Y}(0)\big)\big\rangle\bigg\}.
\end{align*}
The latter expression is not positive thanks to the convexity of the mappings $ (x,y,z,u)\mapsto H(t,x,y,z,u, \hat{p}(t), \hat{q}(t), \hat{k}(t)\big)$ and  $y\mapsto a(y),$
which completes the proof.
 \hfill $\square$

\begin{remark}
In contrast to Agram and {\O}ksendal \cite{Agram-Oksendal} (cf. Theorem 2.1 therein),
Theorem \ref{thm:maxpri} needs not to postulate
the additional joint transverality conditions
\begin{align*}
&\liminf_{T\rightarrow\infty}E\big[e^{\lambda\, T}\langle \hat{Y}(T)-Y(T),\hat{p}(T)\rangle\big]\leq 0, \\[8pt]
& \liminf_{T\rightarrow\infty}E\big[e^{\lambda\, T}\langle\hat{X}(T)-X(T),\hat{q}(T)\rangle\big]\leq 0,
\end{align*}
for every admissible control process $u(\cdot)\in \mathcal{U}_{\lambda},$ since these conditions are implied
by the solvability results of the FBSDEs
(\ref{FBSDEu}) and (\ref{FBSDEpq})  (cf. Theorem \ref{thm:exuniqgeneral} and Remark \ref{rem:exisuniqsigma}), as it is immediately seen in the proof.
\end{remark}

\section{Conclusion}

In this paper we revisit the class of infinite time-horizon FBSDEs given by (\ref{GFBSDEs}),
whose formulation is motivated through several models widely used in economic theory and finance. In particular, we revise the solvability result of Yin \cite{Yin}
by suggesting a more general condition for the parameter $\lambda$ of inequality (\ref{INEQL}). We provide also a comparison result for adapted solutions and establish their continuous dependence on a parameter.
Moreover, we extend the framework of viscosity solutions to the case of random coefficients. Specifically, we introduce the notion of a stochastic viscosity solution for the
quasilinear parabolic BSPDE that is associated with a finite time-horizon system of FBSDEs and ensure its existence by providing a probabilistic representation of it via the
unique solution of this FBSDEs system. Under additional stationary conditions on the random coefficients, we are able to consider next the infinite time-horizon case, in which the stochastic viscosity solution
becomes stationary as well.
Finally, as an application of the solvability of the infinite time-horizon FBSDEs (\ref{GFBSDEs}), we derive a stochastic maximum principle for the stochastic control problem of fully coupled systems of this type.

Last but not least, as far as future research is concerned, this paper motivates highly the investigation of  existence and uniqueness of an either classical or a stochastic viscosity solution for the quasilinear parabolic random
BSPDE (\ref{equ:BSPDEu}), in the sense of Definition \ref{def:svs}. If there exists a unique classical solution, Theorem \ref{thm:class-visc} implies that this solution is the representing random field one may use in order
to decouple and eventually solve the associated system of FBSDEs (\ref{eqn:FBSDEsanss-FINITE}). On the other hand, by the same theorem, in case a unique stochastic viscosity solution exists then it connects the unknown
initial value of the BSDE in terms of the known initial value of the FSDE of the previous system. Therefore, this initial-terminal value system converts to an initial only value system that one may solve by
employing numerical shooting methods.

%\vskip1cm

\end{document}